\documentclass{article}
\newtheorem{Theorem}{Theorem}
\newtheorem{Lemma}[Theorem]{Lemma}
\newtheorem{Definition}[Theorem]{Definition}
\newtheorem{Example}[Theorem]{Example}
\newtheorem{Claim}[Theorem]{Claim}
\begin{document}
\title{Subgraphs of $4$-regular planar graphs}
\author{Chris Dowden \and Louigi Addario-Berry}
\maketitle
\setlength{\unitlength}{1cm}

\begin{abstract}

We shall present an algorithm for determining whether or not a given planar graph $H$ 
can ever be a subgraph of a $4$-regular planar graph.
The algorithm has running time $O \left( |H|^{2.5} \right)$
and can be used to find an explicit $4$-regular planar graph $G \supset H$ if such a graph exists.
It shall not matter whether we specify that $H$ and $G$ must be simple graphs
or allow them to be multigraphs.

\end{abstract}

\section{Introduction}

The last four years have brought large developments
in our understanding of random planar graphs.
Due to the work of a host of people
(see for example~\cite{bod},~\cite{ger} and~\cite{mcd}),
but in particular the landmark paper by Gim\'enez and Noy~\cite{gim},
the generation of random planar graphs and properties such as
the number of edges,
average degree,
and the number of components
are all now well understood.
One such result (from \cite{mcd})
that is related to the subject of this paper is the following:
for any fixed planar graph $H$,
the uniform random planar graph $P_{n}$ on $n$ vertices contains $\Omega (n)$ copies of $H$
except on a set of probability $e^{- \Omega (n)}$.

More recently,
Dowden~\cite{dow} has considered random planar graphs with degree constraints.
For example,
fix $d_{1},d_{2},D_{1}$ and $D_{2}$ and take a graph $P_{n,d_{1},d_{2},D_{1},D_{2}}$
uniformly at random from the set of all planar graphs on $\{ 1,2, \ldots, n \}$
with minimum degree between $d_{1}$ and $d_{2}$, inclusive,
and maximum degree between $D_{1}$ and $D_{2}$, inclusive.
If $D_{2} \geq 3$
(which is necessary for the planarity condition to have any impact),
then it is shown that for any fixed connected planar graph $H$ with $\Delta(H) \leq D_{2}$,
$P_{n,d_{1},d_{2},D_{1},D_{2}}$ also contains $\Omega (n)$ copies of $H$ 
except on a set of probability $e^{- \Omega (n)}$
(where we ignore odd $n$ if $d_{1}=D_{2} \in \{ 3,5 \}$),
apart from two special cases.
The first is when $H$ is $D_{2}$-regular,
in which case the probability is bounded away from both $0$ and $1$,
and the second is when $d_{1}=D_{2}=4$
and $H$ happens to be a graph that can \textit{never} be contained within a $4$-regular planar graph.
This hence raises the question of which graphs can ever be contained in a $4$-regular planar graph
(we will hereafter refer to such graphs as $4$-embeddable),
and that is the topic of this paper. (We mention in passing that there is a related body of work on finding minimal regular 
supergraphs when the planarity restriction is removed; see \cite{bodvan} and the references therein.) 

Note that the problem of determining whether or not a given planar graph $H$ 
can ever be a subgraph of a $k$-regular planar graph is very straightforward for $k \neq 4$,
since the answer is always `yes' if $\Delta (H) \leq k$.
This is clear for $k \in \{ 1,2 \}$,
and can easily be proven for $k \in \{ 3,5 \}$
(note that we must have $k<6$ for planarity)
by showing that in these cases there exist planar graphs 
that are $k$-regular except for exactly one vertex with degree $k-1$,
and that we can hence extend $H$ into a $k$-regular planar graph 
simply by attaching an appropriate number of these graphs to any vertices of $H$ that have degree less than $k$.
This trick does not work for $k=4$, however,
since clearly a graph that is $4$-regular except for exactly one vertex of degree $3$ would have to have
an odd sum of degrees!
In fact,
there do actually exist some planar graphs with maximum degree at most $4$
that are not $4$-embeddable (see Example~\ref{example}), 
and so the matter of determining $4$-embeddability
is non-trivial.

We shall shortly see (in Lemma~\ref{introlemma})
that a given simple planar graph $H$ is actually $4$-embeddable in the world of simple graphs if and only if
it is $4$-embeddable in the world of multigraphs.
Clearly, this second interpretation is just a special case of the more general problem of
determining whether or not a given planar \textit{multigraph} $H$ is $4$-embeddable.
Hence, in this paper we will actually aim to produce an efficient algorithm for the latter problem.

We shall give the details of the algorithm in Section 3.
Before this, in Section 2,
we will prove three lemmas that shall play important roles.
In the first (Lemma~\ref{lemma}),
we shall observe that if we can extend our multigraph $H$ into a $4$-regular planar multigraph,
then we can actually do so without introducing any new vertices.
In the others (Lemmas~\ref{new} and~\ref{lemma2}),
we will show that our problem is straightforward for graphs with a special structure.
In Section 3, we shall then give the algorithm itself,
which will essentially consist of breaking $H$ up
into more and more highly connected pieces until we can apply Lemma~\ref{lemma2}. \\

\begin{figure} [ht] 
\setlength{\unitlength}{1cm}
\begin{picture}(10,2)(-5,0)

\put(0,0){\line(1,0){2}}
\put(1,0.666667){\line(0,1){1.333333}}
\put(0,0){\line(3,2){1}}
\put(0,0){\line(3,4){1}}
\put(0,0){\line(1,2){1}}
\put(2,0){\line(-3,2){1}}
\put(2,0){\line(-3,4){1}}
\put(2,0){\line(-1,2){1}}

\put(0,0){\circle*{0.1}}
\put(2,0){\circle*{0.1}}
\put(1,0.666667){\circle*{0.1}}
\put(1,1.333333){\circle*{0.1}}
\put(1,2){\circle*{0.1}}

\put(-0.1,-0.3){$x$}
\put(2,-0.3){$y$}
\put(0.9,0.366667){$w$}
\put(1.08,1.3){$v$}
\put(1,2.1){$u$}

\end{picture}
\caption{The unique planar embedding of $K_5-uw$.} \label{PH}
\end{figure}
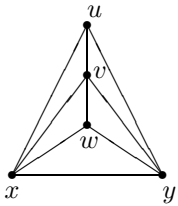

We start with our aforementioned example of a graph with maximum degree~$4$ that is not $4$-embeddable:

\begin{Example} \label{example}
No $4\textrm{-regular}$ planar graph contains a copy of the graph $K_{5}$ minus an edge.
\end{Example}
\textbf{Proof}
The graph $K_{5} - uw$ is drawn with its \textit{unique} planar embedding (see~\cite{whi})
in Figure~\ref{PH}.

Consider any planar graph $G \supset K_{5} -uw$ with $\Delta (G) = 4$.
Since we already have 
$\textrm{deg}_{H}(v) = \textrm{deg}_{H}(x) = \textrm{deg}_{H}(y) = 4$,
any new edge with at least one endpoint inside the triangle given by $vxy$
must have both endpoints inside.
Hence, the sum of degrees inside this triangle must remain odd,
and so this region must still contain a vertex of odd degree.
Thus, $G$ is not $4$-regular.
$\phantom{qwerty}
\setlength{\unitlength}{.25cm}
\begin{picture}(1,1)
\put(0,0){\line(1,0){1}}
\put(0,0){\line(0,1){1}}
\put(1,1){\line(-1,0){1}}
\put(1,1){\line(0,-1){1}}
\end{picture}$ \\

\begin{figure} [ht] 
\setlength{\unitlength}{1cm}
\begin{picture}(10,2)(-5,-1)

\put(-2,0){\line(1,0){6}}

\put(1,1){\line(-1,-1){1}}
\put(1,1){\line(-1,-3){0.333333}}
\put(1,1){\line(1,-3){0.333333}}
\put(1,1){\line(1,-1){1}}

\put(1,-1){\line(-1,1){1}}
\put(1,-1){\line(-1,3){0.333333}}
\put(1,-1){\line(1,3){0.333333}}
\put(1,-1){\line(1,1){1}}

\put(0,0){\circle*{0.1}}
\put(0.666667,0){\circle*{0.1}}
\put(1.333333,0){\circle*{0.1}}
\put(2,0){\circle*{0.1}}
\put(1,1){\circle*{0.1}}
\put(1,-1){\circle*{0.1}}

\put(-2,0){\circle*{0.1}}
\put(4,0){\circle*{0.1}}

\put(-2.1,0.2){$u$}
\put(3.9,0.2){$v$}

\end{picture}
\caption{$\textrm{Constructing a $4$-regular simple planar graph from a $4$-regular planar}$} 
$\textrm{multigraph.}$
\label{4regplanar}
\end{figure}
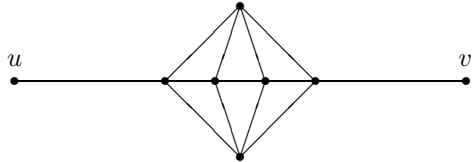

Note that it did not matter whether we took the graph $G$ to be a simple graph or a multigraph.
We shall now conclude this introductory section by seeing that this is indeed always the case:

\begin{Lemma} \label{introlemma}
Given a simple planar graph $H$,
there exists a $4$-regular simple planar graph $G \supset H$
if and only if there exists a $4$-regular planar multigraph $G^{\prime} \supset H$.
\end{Lemma}
\textbf{Proof}
The forward implication is trivial.
To see the reverse implication,
replace every edge $e=uv$ of $E(G^{\prime}) \setminus E(H)$ by a copy of the graph shown in Figure~\ref{4regplanar}.
\\
The resulting graph $G$ will be a $4$-regular simple planar graph with $H \subset G$.
$\phantom{w}
\setlength{\unitlength}{.25cm}
\begin{picture}(1,1)
\put(0,0){\line(1,0){1}}
\put(0,0){\line(0,1){1}}
\put(1,1){\line(-1,0){1}}
\put(1,1){\line(0,-1){1}}
\end{picture}$ \\
\\
\\

\section{Lemmas and Definitions}

In this section, we shall do the groundwork for our algorithm by proving three easy but important lemmas,
as well as introducing several helpful definitions. For $v \in V(G)$, $e \in E(G)$, we write $G-v$ for the graph $(V(G)\setminus\{v\},\{e \in E(G): v \not\in e\})$, and define $G+v, G-e$, and $G+e$ similarly. 

We start by showing that if there exists a $4$-regular planar multigraph~$G \supset~H$,
then we may assume that $V(G)=V(H)$:

\begin{Lemma} \label{lemma}
Given a planar multigraph $H$,
there exists a $4$-regular planar multigraph $G \supset H$
if and only if there exists a $4$-regular planar multigraph $G^{\prime} \supset H$ 
with $V(G^{\prime}) = V(H)$.
\end{Lemma}
\textbf{Proof}
The `if' direction is trivial,
so it will suffice to show the `only if' direction.

Suppose there exists a $4$-regular planar multigraph $G \supset H$,
and let $G^{\prime}$ be a \textit{minimal} such graph,
in the sense that $|V(G^{\prime}) \setminus V(H)|$ is as small as possible.

Suppose $|V(G^{\prime}) \setminus V(H)| \neq 0$
(hoping to obtain a contradiction)
and let $v \in~V(G^{\prime}) \setminus V(H)$.
We shall show that we can obtain a $4$-regular planar multigraph $G^{*}$
such that $H \subset G^{*}$ and $V(G^{*}) = V(G^{\prime}) \setminus v$,
thus obtaining our desired contradiction: \\
Case (a): If $v$ has two loops to itself,
then we may simply take $G^{*}$ to be $G^{\prime} -v$. \\
Case (b): If $v$ has exactly one loop to itself
and its other two neighbours are $v_{1}$ and $v_{2}$
(where we allow the possibility that $v_{1}=v_{2}$),
then we may take $G^{*}$ to be $(G^{\prime} - v)+ v_{1}v_{2}$. \\
Case (c): If $v$ has no loops to itself,
then fix a plane drawing of $G^{\prime}$
and let $e_{1},e_{2},e_{3}$ and $e_{4}$ be the edges incident to $v$ 
\textit{in clockwise order} in this plane drawing.
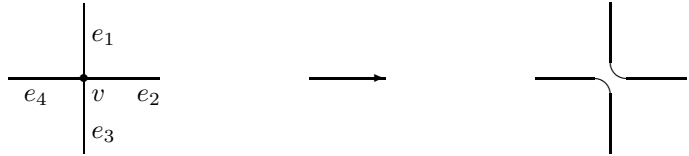
\begin{figure} [ht] 
\setlength{\unitlength}{1cm}
\begin{picture}(10,2)(-1.5,0)

\put(0,1){\line(1,0){2}}
\put(1,0){\line(0,1){2}}
\put(1,1){\circle*{0.1}}
\put(1.1,1.5){$e_{1}$}
\put(1.7,0.7){$e_{2}$}
\put(1.1,0.2){$e_{3}$}
\put(0.2,0.7){$e_{4}$}
\put(1.1,0.7){$v$}

\put(4,1){\vector(1,0){1}}

\put(7,1){\line(1,0){0.75}}
\put(8.25,1){\line(1,0){0.75}}
\put(8,0){\line(0,1){0.75}}
\put(8,1.25){\line(0,1){0.75}}
\put(8.25,1.25){\oval(0.5,0.5)[bl]}
\put(7.75,0.75){\oval(0.5,0.5)[tr]}

\end{picture}
\caption{$\textrm{Constructing a smaller $4$-regular planar multigraph in case (c).}$} \label{plane}
\end{figure} 
Let $v_{1},v_{2},v_{3}$ and $v_{4}$, respectively, denote the other endpoints of $e_{1},e_{2},e_{3}$ and $e_{4}$
(allowing the possibility that $v_{i}=v_{j}$ for some $i$ and $j$).
Then we may take $G^{*}$
to be $(G^{\prime} -v)+ v_{1}v_{2} + v_{3}v_{4}$,
since this can also be drawn in the plane (see Figure~\ref{plane}).
$\phantom{qwerty}
\setlength{\unitlength}{.25cm}
\begin{picture}(1,1)
\put(0,0){\line(1,0){1}}
\put(0,0){\line(0,1){1}}
\put(1,1){\line(-1,0){1}}
\put(1,1){\line(0,-1){1}}
\end{picture}$ \\

Note that Lemma~\ref{lemma} itself provides a way to determine algorithmically
whether or not a given planar multigraph $H$ can ever be a subgraph of a $4$-regular planar multigraph,
since it will suffice just to check all $4$-regular planar multigraphs with the same vertex set as $H$.
However, we shall produce a much faster algorithm in Section 3. 
The following lemma is a step in that direction,
establishing a polynomial time algorithm for the case when $H$ is $3$-vertex-connected:

\begin{Lemma} \label{new}
Given any $3$-vertex-connected planar multigraph $H$,
we can determine in $O \left( |H|^{2.5} \right)$ operations whether or not $H$ is $4$-embeddable.
\end{Lemma}
\textbf{Proof}
Without loss of generality,
$\Delta(H) \leq 4$.
Thus, $3$-vertex-connectivity implies that $H$ has no loops.
By a result of Whitney~\cite{whi} on $3$-vertex-connected simple graphs,
it then follows that $H$ has a unique planar embedding
\footnote{By ``unique'' we mean ``unique up to automorphisms of plane graphs''. In particular, the presence of parallel edges cannot affect the uniqueness of the embedding.},
and this can be obtained in $O(|H|)$ operations (see~\cite{boo}).
We shall use this embedding to reduce the problem of $4$-embeddability to finding a perfect matching
in a suitably defined `auxiliary' graph.

Let $f(v) = 4 -\deg_{H}(v)$~$\forall v \in V(H)$.
Then we define our auxiliary graph $A$
(which will not necessarily be planar)
to consist of $f(v)$ copies, $v_{1},v_{2}, \ldots, v_{f(v)}$,
of each vertex $v$,
with an edge between $v_{i}$ and $w_{j}$ 
precisely if $v$ and $w$ are on a common face in $H$
(here we allow $v=w$, but not $v_{i}=w_{j}$).
We claim that $H$ is $4$-embeddable if and only if $A$ has a perfect matching.

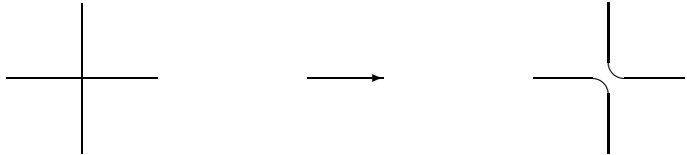
\begin{figure} [ht] 
\setlength{\unitlength}{1cm}
\begin{picture}(10,2)(-1.5,0)

\put(0,1){\line(1,0){2}}
\put(1,0){\line(0,1){2}}

\put(4,1){\vector(1,0){1}}

\put(7,1){\line(1,0){0.75}}
\put(8.25,1){\line(1,0){0.75}}
\put(8,0){\line(0,1){0.75}}
\put(8,1.25){\line(0,1){0.75}}
\put(8.25,1.25){\oval(0.5,0.5)[bl]}
\put(7.75,0.75){\oval(0.5,0.5)[tr]}

\end{picture}
\caption{$\textrm{Separating crossing edges of our matching.}$}
\label{D15}
\end{figure}

First,
suppose $H$ is $4$-embeddable,
i.e.~$H \subset G$ for some $4$-regular planar multigraph $G$.
By Lemma~\ref{lemma},
we may assume that $V(G)=V(H)$,
in which case the edges $E(G) \setminus E(H)$
form a perfect matching in $A$
(where we choose which copy of $v$ to use for a given edge in some arbitrary consistent manner).
Conversely, if we can find a perfect matching in $A$,
then inserting the edges of this matching into our embedding of $H$ will give us a
(not necessarily plane) $4$-regular multigraph,
which can then be made into a plane $4$-regular multigraph simply by separating any new edges that cross,
as in Figure~\ref{D15}.

Note that we can obtain the auxiliary graph in $O \left( |H|^{2} \right)$ time,
since there are $O \left( |H|^{2} \right)$ possible edges.
This graph will have at most $4|H|$ vertices,
and so we can then determine whether or not it has a perfect matching in $O \left( |H|^{2.5} \right)$ time (see~\cite{eve}).
$\phantom{qwerty}
\setlength{\unitlength}{.25cm}
\begin{picture}(1,1)
\put(0,0){\line(1,0){1}}
\put(0,0){\line(0,1){1}}
\put(1,1){\line(-1,0){1}}
\put(1,1){\line(0,-1){1}}
\end{picture}$ \\

We remark that the only place in the preceding proof where we used $3$-vertex-connectivity
was to ensure that $H$ had a unique planar embedding.
Thus, an identical proof gives a polynomially testable necessary and sufficient condition
for any given \textit{plane} graph to be a subgraph of a $4$-regular planar graph.
However, note that this condition does not yield an efficient algorithm for an arbitrary \textit{planar} graph $H$,
as in general $H$ may have exponentially many planar embeddings. \\

Returning to our main thrust,
recall that in the proof of Lemma~\ref{new} we introduced the function $f$,
which encoded the discrepancy between a vertex's degree in $H$ and its target degree
(in this case, there was a target degree of $4$ for all vertices).
Such functions will play an important role in our algorithm,
and so we will now take the time to set up a more general framework for them:

\begin{Definition}
Given a planar multigraph $H$ with maximum degree at most $4$,
$f_{H}:V(H) \to \mathbf{N}$
is a \emph{discrepancy function} on $H$
if (a) $f_{H}(v) \leq 4-\deg_{H}(v)$~$\forall v \in V(H)$
(we call this \emph{the discrepancy inequality})
and (b) $\sum_{v \in V(H)} f_{H}(v)$ is even
(we call this \emph{discrepancy parity}).
If it is also the case that $f_{H}(v) + \deg_{H}(v)$ is even for all $v \in V(H)$,
we call $f_{H}$ an \emph{even discrepancy function} on $H$.

We say that a plane multigraph $G$ \emph{satisfies} $(H,f_{H})$
if $V(G)=V(H), E(G) \supset E(H)$ and $\deg_{G}(v) = \deg_{H}(v) + f_{H}(v)$~$\forall v$.
If such a plane multigraph $G$ exists,
we say that \emph{$f_{H}$ can be satisfied on $H$},
or that \emph{$(H,f_{H})$ can be satisfied}. 
\end{Definition}

We next also introduce `augmentations',
a graph operation we will use repeatedly during the algorithm:

\begin{Definition}
Given a multigraph $B$,
we define the operation of \emph{placing a diamond}
on an edge $uv \in E(B)$
to mean that we subdivide the edge with three vertices 
and then also add two other new vertices so that they are both adjacent to precisely these three vertices.
We define the operation of \emph{placing a vertex}
on an edge $xy \in E(B)$
to mean that we subdivide the edge with a single vertex.

Given multigraphs $B$ and $R$
and a discrepancy function $f_{R}$,
we say that $(R,f_{R})$ is an \emph{augmentation} of $B$ if $R$ can be formed from $B$ 
by placing vertices and diamonds on some of the edges of $B$
(in such a way that there is at most one vertex or diamond on each original edge)
and if $f_{R} = 4 - \deg_{R}$ for all vertices in the new diamonds
and $f_{R} \in \{1,2\}$ for the other new vertices.
\end{Definition}

\begin{figure} [ht] 
\setlength{\unitlength}{1cm}
\begin{picture}(10,2)(-1.85,0)

\put(0,0){\line(1,0){2}}
\put(0,2){\line(1,0){2}}
\put(0,1){\oval(1,2)}

\put(6,0){\line(1,0){2}}
\put(6,2){\line(1,0){2}}
\put(6,1){\oval(1,2)}
\put(5,1){\line(1,1){0.5}}
\put(5,1){\line(1,-1){0.5}}
\put(5,1){\line(1,0){1}}
\put(6,1){\line(-1,1){0.5}}
\put(6,1){\line(-1,-1){0.5}}

\put(0,0){\circle*{0.1}}
\put(0,2){\circle*{0.1}}
\put(2,0){\circle*{0.1}}
\put(2,2){\circle*{0.1}}

\put(6,0){\circle*{0.1}}
\put(6,2){\circle*{0.1}}
\put(7,0){\circle*{0.1}}
\put(8,0){\circle*{0.1}}
\put(8,2){\circle*{0.1}}
\put(6.5,1){\circle*{0.1}}
\put(5,1){\circle*{0.1}}
\put(5.5,1){\circle*{0.1}}
\put(6,1){\circle*{0.1}}
\put(5.5,0.5){\circle*{0.1}}
\put(5.5,1.5){\circle*{0.1}}

\put(2.3,1.75){\large{$B$}}
\put(8.3,1.75){\large{$(R,f_{R})$}}

\put(5.9,-0.4){$0$}
\put(5.9,2.2){$1$}
\put(6.9,-0.4){$1$}
\put(7.9,-0.4){$3$}
\put(7.9,2.2){$1$}
\put(6.6,0.9){$2$}
\put(4.7,0.9){$1$}
\put(5.55,1){$0$}
\put(6.1,0.9){$1$}
\put(5.3,0.2){$0$}
\put(5.3,1.6){$0$}

\end{picture}
\caption{$\textrm{A planar multigraph and an augmentation of it.}$}
\label{D1}
\end{figure}
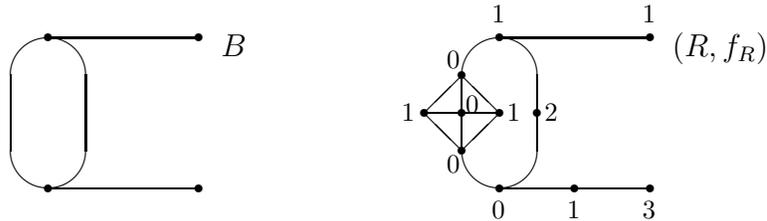

An example of an augmentation is given in Figure~\ref{D1}.
When we break $H$ up into pieces in our algorithm,
the augmentation of a piece will capture the key information about how it interacted with the rest of $H$. \\

Note that in performing such an augmentation, we may differently augment two parallel edges (as in Figure \ref{D1}). This causes an asymmetry in the graph, and so $R$ 
may not have a unique planar embedding even if $B$ does. However, augmenting two parallel edges in the {\em same} way will not affect the uniqueness of the embedding. 
We will use this fact in Lemma~\ref{lemma2}, below, 
which extends Lemma~\ref{new}
to cover augmentations of $3$-vertex-connected multigraphs:

\begin{Lemma} \label{lemma2}
Let $B$ be a planar multigraph of maximum degree at most $4$ that contains no $2$-vertex-cuts
(so $|B| \leq 3$ or $B$ is $3$-vertex-connected),
and let $(R,f_{R})$ be an augmentation of $B$.
Suppose we know which parts of~$R$ correspond to which edges of $B$.
Then $\exists \lambda$, independent of $B$ and $R$,
such that we can determine in at most $\lambda |B|^{2.5}$ operations
whether or not $(R,f_{R})$ can be satisfied.
\end{Lemma}
\textbf{Proof}
Without loss of generality,
we may assume that $|B|>3$
(since if $B$ is bounded then there are only a finite number of possibilities for $(R,f_{R})$,
and the satisfiability of these can be determined in finite time).
Thus, as in Lemma~\ref{new},
$B$ has no loops and so has a unique planar embedding.
Hence, $R$ will also have a unique planar embedding,
apart from possibly at places where $B$ has multi-edges. 

Note that all vertices in $B$ must have at least $3$ distinct neighbours,
since $B$ does not contain any $2\textrm{-vertex-cuts}$.
Hence (since $\deg_{B}(x) = \deg_{R}(x) \leq 4 - f_{R}(x)$~$\forall x \in V(B)$),
if vertices $u$ and $v$ have a multi-edge between them in~$B$,
then it must be only a double-edge
and it must be that $f_{R}(u) = f_{R}(v) = 0$.
We shall now use this information to find a pair $(R^{\prime},f_{R^{\prime}})$
such that $R^{\prime}$ has a unique planar embedding and
$(R,f_{R})$ can be satisfied if and only if
we can satisfy $(R^{\prime},f_{R^{\prime}})$.

\begin{figure} [ht] 
\setlength{\unitlength}{1cm}
\begin{picture}(10,2)(-3.5,0)

\put(0,0){\line(0,1){2}}
\put(0.1,0){$A$}

\put(1.5,0){\line(0,1){2}}
\put(1.5,1){\circle*{0.1}}
\put(1.6,0.9){$1$}
\put(1.6,0){$B$}

\put(3,0){\line(0,1){2}}
\put(3,1){\circle*{0.1}}
\put(3.1,0.9){$2$}
\put(3.1,0){$C$}

\put(4.5,0){\line(0,1){2}}
\put(4,1){\line(1,0){1}}
\put(4,1){\line(1,1){0.5}}
\put(4,1){\line(1,-1){0.5}}
\put(5,1){\line(-1,1){0.5}}
\put(5,1){\line(-1,-1){0.5}}

\put(4,1){\circle*{0.1}}
\put(4.5,1){\circle*{0.1}}
\put(5,1){\circle*{0.1}}
\put(4.5,0.5){\circle*{0.1}}
\put(4.5,1.5){\circle*{0.1}}

\put(3.7,0.9){$1$}
\put(4.55,1){$0$}
\put(5.1,0.9){$1$}
\put(4.3,0.2){$0$}
\put(4.3,1.6){$0$}

\put(5.1,0){$D$}

\end{picture}
\caption{$\textrm{Augmented versions of an edge.}$}
\label{D2}
\end{figure}
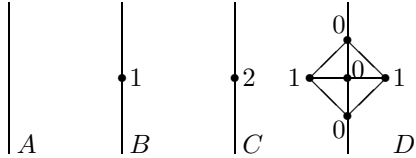

\begin{figure} [ht] 
\setlength{\unitlength}{1cm}
\begin{picture}(10,2.8)(-0.25,-0.4)

\put(6,1){\oval(1,2)}
\put(5,1){\line(1,1){0.5}}
\put(5,1){\line(1,-1){0.5}}
\put(5,1){\line(1,0){1}}
\put(6,1){\line(-1,1){0.5}}
\put(6,1){\line(-1,-1){0.5}}
\put(6,2){\line(1,1){0.4}}
\put(6,2){\line(-1,1){0.4}}
\put(6,0){\line(1,-1){0.4}}
\put(6,0){\line(-1,-1){0.4}}

\put(5,1){\circle*{0.1}}
\put(5.5,1){\circle*{0.1}}
\put(6,1){\circle*{0.1}}
\put(5.5,0.5){\circle*{0.1}}
\put(5.5,1.5){\circle*{0.1}}
\put(6,0){\circle*{0.1}}
\put(6,2){\circle*{0.1}}

\put(4.7,0.9){$1$}
\put(5.55,1){$0$}
\put(6.1,0.9){$1$}
\put(5.3,0.2){$0$}
\put(5.3,1.6){$0$}
\put(5.92,2.15){$0$}
\put(5.92,-0.35){$0$}

\put(6.2,2){$u$}
\put(6.2,-0.15){$v$}

\end{picture}
\caption{$\textrm{A Type A-Type D double edge.}$}
\label{D3}
\end{figure}
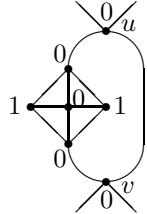

Let Type A, Type B, Type C and Type D denote the four possible `augmented versions' of an edge,
as shown in Figure~\ref{D2},
and recall that $R$ will have a unique embedding apart from at any places where $B$ has a double-edge. 
More strongly, the only possible non-uniqueness arises when $B$ has a double-edge and these edges are of distinct types
in the augmentation. (This is explained in the paragraph just before the statement of the lemma.) We now explain how to deal with such situations. 
If there exist vertices $u$ and $v$ with a Type A-Type D double-edge between them,
then it can be seen that it is impossible to satisfy $(R,f_{R})$,
since $f(u)=f(v)=0$ (see Figure~\ref{D3}).
If we have no Type A-Type D double-edges,
then let $R^{\prime}$ be formed from $R$ as follows: \\
(i) If the augmented versions of a double-edge are Type A and Type B,
then delete the Type A part; \\
(ii) If A and C, delete C; \\
(iii) If B and C, delete C; \\
(iv) If B and D, delete D. \\
(v) If C and D, delete C. \\
Let $f_{R^{\prime}}(v) = f_{R}(v)$~$\forall v \in V(R^{\prime})$.

Using the fact that the two ends of any double-edge must have $f_{R} = 0$,
it is easy to see that $(R,f_{R})$ can be satisfied if and only if $(R^{\prime},f_{R^{\prime}})$ can be satisfied.
It is also clear that $R^{\prime}$ will have a unique embedding
(if we consider the vertices as unlabelled apart from the discrepancy function).
Thus, to determine whether or not $(R^{\prime},f_{R^{\prime}})$ can be satisfied,
it suffices to see if we can satisfy $(R^{\prime},f_{R^{\prime}})$ in this embedding.
 
Note that $R^{\prime}$ can be generated from $R$ in $O \left( |B| \right)$ time,
since there are $O \left( |B| \right)$ edges in $B$
and we know which parts of $R$ correspond to which edges of $B$.
Since $|R^{\prime}|=O(|B|)$,
the planar embedding of $R^{\prime}$ can then be found in $O(|B|)$ time (see~\cite{boo}),
and the satisfiability of $(R^{\prime},f_{R^{\prime}})$ can then be determined in $O \left( |B|^{2.5} \right)$
operations by looking for a perfect matching in the corresponding auxiliary graph,
as in the proof of Lemma~\ref{new}.
$\phantom{q}
\setlength{\unitlength}{.25cm}
\begin{picture}(1,1)
\put(0,0){\line(1,0){1}}
\put(0,0){\line(0,1){1}}
\put(1,1){\line(-1,0){1}}
\put(1,1){\line(0,-1){1}}
\end{picture}$ \\

It is worth remarking that,
apart from obtaining the auxiliary graph and looking for a perfect matching,
all the procedures involved in the previous proof actually only take $O(|B|)$ time. \\
\\
\\

\section{The Algorithm}

We will now present our algorithm. 
We shall first provide a short sketch,
before then giving the details in full.
Afterwards,
we will investigate the running time.

\subsection*{Sketch of Algorithm} 
The algorithm shall consist of four stages,
each of which will involve breaking $H$ up into more highly connected pieces,
until we can eventually apply Lemma~\ref{lemma2} to all of these.

We will start,
in Stages 1 and 2,
by straightforwardly showing that $H$ is $4$-embeddable
if and only if all its $2$-edge-connected components are.
This part of our argument will not require the use of either Lemma~\ref{lemma} or Lemma~\ref{lemma2}.

In Stage 3, 
we will break our $2$-edge-connected components into $2$-\textit{vertex}-connected blocks,
and show that the discrepancy function $f=4-\deg$ can be satisfied on our $2$-edge-connected components
if and only if certain specified discrepancy functions can be satisfied on all the $2$-vertex-connected blocks.
We shall use Lemma~\ref{lemma} to simplify our arguments here.

Stage 4 is where we will use the notion of augmentations.
We shall split our $2$-vertex-connected blocks into $3$-vertex-connected multigraphs
and define augmented versions of each of these.
There will be different cases depending on exactly how the $2$-vertex-cuts break up the graph,
and we will show that the discrepancy functions defined on our $2$-vertex-connected blocks can be satisfied
if and only if all these augmentations can be satisfied.
This can then be determined using Lemma~\ref{lemma2}.

\subsection*{FULL ALGORITHM}

\subsubsection*{STAGE 1}

Clearly, there exists a $4$-regular planar multigraph $G \supset H$ if and only if
there exist $4$-regular planar multigraphs $G_{i} \supset H_{i}$ for all components $H_{i}$ of $H$
(the `if' direction follows by taking $G$ to be the graph whose components are the $G_{i}$'s 
and the `only if' direction follows by taking $G_{i}=G$~$\forall i$).

Thus, the first stage of our algorithm will be to split $H$ into its components.

\subsubsection*{STAGE 2}

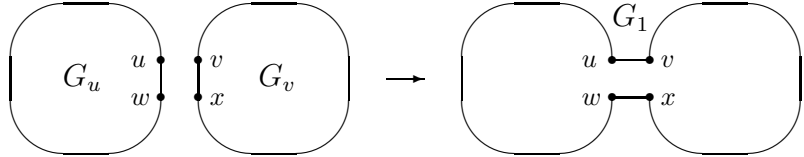
\begin{figure} [ht] 
\setlength{\unitlength}{1cm}
\begin{picture}(10,2)(-0.75,0)

\put(1,1){\oval(2,2)}
\put(3.5,1){\oval(2,2)}
\put(7,1.25){\oval(2,1.5)[t]}
\put(7,0.75){\oval(2,1.5)[b]}
\put(9.5,1.25){\oval(2,1.5)[t]}
\put(9.5,0.75){\oval(2,1.5)[b]}

\put(5,1){\vector(1,0){0.5}}

\put(8,0.75){\line(1,0){0.5}}
\put(8,1.25){\line(1,0){0.5}}
\put(6,0.75){\line(0,1){0.5}}
\put(10.5,0.75){\line(0,1){0.5}}

\put(2,0.75){\circle*{0.1}}
\put(2,1.25){\circle*{0.1}}
\put(2.5,0.75){\circle*{0.1}}
\put(2.5,1.25){\circle*{0.1}}
\put(8,0.75){\circle*{0.1}}
\put(8,1.25){\circle*{0.1}}
\put(8.5,0.75){\circle*{0.1}}
\put(8.5,1.25){\circle*{0.1}}

\put(1.6,1.15){$u$}
\put(1.6,0.65){$w$}
\put(2.65,1.15){$v$}
\put(2.65,0.65){$x$}

\put(7.6,1.15){$u$}
\put(7.6,0.65){$w$}
\put(8.65,1.15){$v$}
\put(8.65,0.65){$x$}

\put(0.7,0.9){\large{$G_{u}$}}
\put(3.3,0.9){\large{$G_{v}$}}
\put(8,1.7){\large{$G_{1}$}}

\end{picture}
\caption{$\textrm{Constructing a $4$-regular planar multigraph $G_{1}$ 
from $4$-regular planar}$}
\textrm{multigraphs $G_{u}$ and $G_{v}$.}
\label{D4}
\end{figure}

Let $H_{1}$ be a component of $H$ and suppose that $H_{1}$ has a cut-edge $e=uv$. 
Let $H_{u}$ and $H_{v}$ denote the components of $H_{1} - e$ containing $u$ and $v$, respectively.
Clearly, there exists a $4$-regular planar multigraph $G_{1} \supset H_{1}$ only if
there exist $4$-regular planar multigraphs $G_{u} \supset H_{u}$ and $G_{v} \supset H_{v}$
(this follows by taking $G_{u}=G_{v}=G_{1}$).
We shall now see that the converse is also true:

Suppose there exist $4$-regular planar multigraphs $G_{u} \supset H_{u}$ and $G_{v} \supset H_{v}$.
Note that $\deg_{H_{u}}(u) = \deg_{H_{1}}(u)-1 \leq 3$,
since $v \notin V(H_{u})$,
so $\exists w \in V(G_{u})$ such that
$uw \in E(G_{u}) \setminus E(H_{u})$. (It is possible that $w=u$.) 
Similarly,
$\exists x \in V(G_{v})$ such that
$vx \in E(G_{v}) \setminus E(H_{v})$.
Since $G_{u}$ and $G_{v}$ are both planar,
they can be drawn with the edges $uw$ and $vx$, respectively, in the outside face.
Thus, the graph $G_{1}$ formed by deleting these two edges and inserting edges $uv$ and $wx$ will also be planar,
as well as being a $4$-regular multigraph containing $H_{1}$ (see Figure~\ref{D4}).

We have shown that $H_{1}$ is $4$-embeddable if and only if $H_{u}$ and $H_{v}$ both are.
Thus, by repeated use of this result,
we find that $H_{1}$ is $4$-embeddable if and only if
all its $2$-edge-connected components are
(counting an isolated vertex as $2$-edge-connected).

Therefore, the second stage of our algorithm will be to split the components of $H$ 
into their $2$-edge-connected components.

\subsubsection*{STAGE 3}

Let $A$ be one of our $2$-edge-connected components.
We wish to determine whether or not there exists a $4$-regular planar multigraph $G_{A} \supset A$.
By Lemma~\ref{lemma},
it suffices to discover whether or not there exists a $4$-regular planar multigraph $G^{\prime} \supset A$
with $V(G^{\prime}) = V(A)$,
i.e.~to determine whether or not we can satisfy the even discrepancy function on $A$ defined by setting
$f_{A}(v) = 4-\deg_{A}(v)$~$\forall v \in V(A)$.

Suppose that $A$ has a cut-vertex $v$.
Since $A$ contains no cut-edges,
it must be that $A - v$ consists of exactly two components, $A_{1}$ and $A_{2}$,
with exactly two edges from $v$ to each of these components.
Thus, $\deg_{A}(v)=4$ and $f_{A}(v)=0$.

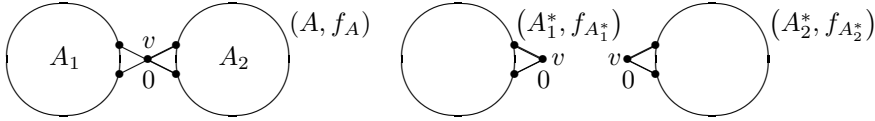
\begin{figure} [ht] 
\setlength{\unitlength}{0.75cm}
\begin{picture}(10,2)(-0.3,0)

\put(1,1){\oval(2,2)}
\put(4,1){\oval(2,2)}
\put(8,1){\oval(2,2)}
\put(12.5,1){\oval(2,2)}

\put(2,0.75){\line(2,1){1}}
\put(2,1.25){\line(2,-1){1}}
\put(9,0.75){\line(2,1){0.5}}
\put(9,1.25){\line(2,-1){0.5}}
\put(11,1){\line(2,1){0.5}}
\put(11,1){\line(2,-1){0.5}}

\put(5,1.5){$(A,f_{A})$}
\put(9,1.5){$\left( A_{1}^{*},f_{A_{1}^{*}} \right)$}
\put(13.5,1.5){$\left( A_{2}^{*},f_{A_{2}^{*}} \right)$}

\put(2.5,1){\circle*{0.1333}}
\put(9.5,1){\circle*{0.1333}}
\put(11,1){\circle*{0.1333}}

\put(2,0.75){\circle*{0.1333}}
\put(2,1.25){\circle*{0.1333}}
\put(3,0.75){\circle*{0.1333}}
\put(3,1.25){\circle*{0.1333}}

\put(9,0.75){\circle*{0.1333}}
\put(9,1.25){\circle*{0.1333}}
\put(11.5,0.75){\circle*{0.1333}}
\put(11.5,1.25){\circle*{0.1333}}

\put(2.4,0.5){$0$}
\put(2.4,1.2){$v$}
\put(9.4,0.5){$0$}
\put(9.65,0.9){$v$}
\put(10.9,0.5){$0$}
\put(10.65,0.9){$v$}

\put(0.75,0.9){$A_{1}$}
\put(3.75,0.9){$A_{2}$}

\end{picture}
\caption{$\textrm{The planar multigraphs $A,A_{1}^{*}$ and $A_{2}^{*}$.}$}
\label{D5}
\end{figure}

Let $A_{1}^{*}$ denote the planar multigraph induced by $V(A_{1}) \cup v$ 
and let $f_{A_{1}^{*}}$ denote the even discrepancy function on $A_{1}^{*}$ defined by setting 
$f_{A_{1}^{*}}(x) = f_{A}(x)$~$\forall x \in V(A_{1}^{*})$
(we have $f_{A_{1}^{*}}(x) + \deg_{A_{1}^{*}}(x) = 4$~$\forall x \neq v$
and $f_{A_{1}^{*}}(v) + \deg_{A_{1}^{*}}(v) = 2$,
so $f_{A_{1}^{*}}$ is indeed an even discrepancy function).
Let $A_{2}^{*}$ and $f_{A_{2}^{*}}$ be defined similarly (see Figure~\ref{D5}).

\begin{figure} [ht] 
\setlength{\unitlength}{1cm}
\begin{picture}(10,3.5)(-0.5,-1)

\put(0.75,0.75){\oval(1.5,1.5)}
\put(3.25,0.75){\oval(1.5,1.5)}

\put(1.5,0.5){\line(2,1){1}}
\put(1.5,1){\line(2,-1){1}}

\put(4.2,2.25){\large{$G^{\prime}$}}

\put(2,0.75){\circle*{0.1}}

\put(1.9,0.45){$v$}

\put(1,1.5){\line(1,0){2}}
\put(1,0){\line(1,0){2}}
\put(2,1.5){\oval(3,1)[t]}
\put(2,0.75){\oval(4,3.5)}
\put(2,0){\oval(3,1)[b]}

\put(0.6,0.6){\large{$A_{1}$}}
\put(3.1,0.6){\large{$A_{2}$}}

\put(0.5,1.5){\circle*{0.1}}
\put(1,1.5){\circle*{0.1}}
\put(0,0.5){\circle*{0.1}}
\put(0,1){\circle*{0.1}}
\put(0.5,0){\circle*{0.1}}
\put(1,0){\circle*{0.1}}
\put(1.5,0.5){\circle*{0.1}}
\put(1.5,1){\circle*{0.1}}

\put(3,1.5){\circle*{0.1}}
\put(3.5,1.5){\circle*{0.1}}
\put(2.5,0.5){\circle*{0.1}}
\put(2.5,1){\circle*{0.1}}
\put(3,0){\circle*{0.1}}
\put(3.5,0){\circle*{0.1}}
\put(4,0.5){\circle*{0.1}}
\put(4,1){\circle*{0.1}}

\put(5,0.75){\vector(1,0){1}}

\put(7.75,0.75){\oval(1.5,1.5)}
\put(10.25,0.75){\oval(1.5,1.5)}

\put(8.5,0.5){\line(2,1){1}}
\put(8.5,1){\line(2,-1){1}}

\put(11.2,1.5){\large{$G^{*}$}}

\put(9,0.75){\circle*{0.1}}

\put(7.5,1.5){\circle*{0.1}}
\put(8,1.5){\circle*{0.1}}
\put(7,0.5){\circle*{0.1}}
\put(7,1){\circle*{0.1}}
\put(7.5,0){\circle*{0.1}}
\put(8,0){\circle*{0.1}}
\put(8.5,0.5){\circle*{0.1}}
\put(8.5,1){\circle*{0.1}}

\put(10,1.5){\circle*{0.1}}
\put(10.5,1.5){\circle*{0.1}}
\put(9.5,0.5){\circle*{0.1}}
\put(9.5,1){\circle*{0.1}}
\put(10,0){\circle*{0.1}}
\put(10.5,0){\circle*{0.1}}
\put(11,0.5){\circle*{0.1}}
\put(11,1){\circle*{0.1}}

\put(8.9,0.45){$v$}

\put(7.75,1.5){\oval(0.5,0.5)[t]}
\put(7,0.75){\oval(0.5,0.5)[l]}
\put(7.75,0){\oval(0.5,0.5)[b]}

\put(10.25,1.5){\oval(0.5,0.5)[t]}
\put(11,0.75){\oval(0.5,0.5)[r]}
\put(10.25,0){\oval(0.5,0.5)[b]}

\end{picture}
\caption{$\textrm{Constructing the graph $G^{*}$ from $G^{\prime}$.}$}
\label{D6}
\end{figure}
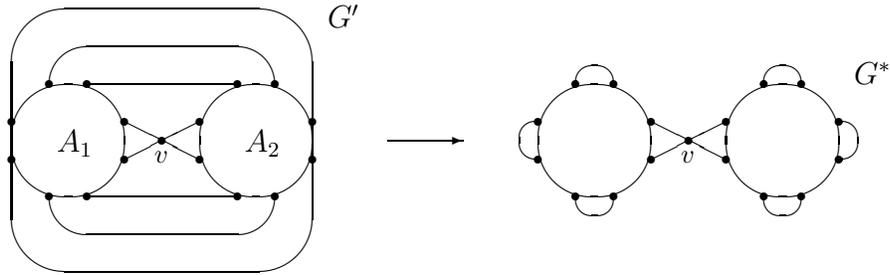

Clearly, we can satisfy $(A,f_{A})$ if we can satisfy both 
$\left( A_{1}^{*},f_{A_{1}^{*}} \right)$ and $\left( A_{2}^{*},f_{A_{2}^{*}} \right)$
(since if there exist plane multigraphs $G_{1}^{*}$ and $G_{2}^{*}$
satisfying $\left( A_{1}^{*},f_{A_{1}^{*}} \right)$ and $\left( A_{2}^{*},f_{A_{2}^{*}} \right)$, respectively,
then we may assume that $v$ is in the outside face of both of these,
and so we can then `glue' these two drawings together at $v$ to obtain a plane multigraph that satisfies $(A,f_{A})$).
We shall now see that the converse is also true:

Suppose $(A,f_{A})$ can be satisfied,
i.e.~there exists a plane multigraph $G^{\prime} \supset A$
with $V(G^{\prime})=V(A)$ and $\deg_{G^{\prime}}(x) = 4$~$\forall x$.
Let us consider the induced plane drawing of $A$.
Since $A_{2}$ is connected,
it must lie in a single face of $A_{1}^{*}$.
Thus, we may assume that our plane drawing of $A$ is as shown in Figure~\ref{D6},
where without loss of generality we have drawn $A_{2}$ in the outside face of $A_{1}^{*}$.
Note that the set of edges in $E(G^{\prime}) \setminus E(A)$
between $A_{1}$ and $A_{2}$ must all lie in a single face of our plane drawing
and that there must be an even number of such edges,
since $f_{A}$ is an even discrepancy function and $f_{A}(v)=0$.
Thus, we may `pair up' these edges,
as in Figure~\ref{D6},
to obtain a plane multigraph $G^{*}$ satisfying $(f_{A},A)$ that has \textit{no} edges from $A_{1}$ to $A_{2}$.
It is then clear that $G_{1}^{*} = G^{*} \setminus A_{2}$ and $G_{2}^{*} = G^{*} \setminus A_{1}$
will satisfy $\left( A_{1}^{*},f_{A_{1}^{*}} \right)$ and $\left( A_{2}^{*},f_{A_{2}^{*}} \right)$, respectively.

Thus, we have shown that the even discrepancy function $f_{A}$ can be satisfied on $A$ if and only if
the even discrepancy functions $f_{A_{1}^{*}}$ and $f_{A_{2}^{*}}$ 
can be satisfied on $A_{1}^{*}$ and $A_{2}^{*}$, respectively.
By repeatedly using this result, 
we may obtain a set of discrepancy functions 
defined on $2$-\textit{vertex}-connected planar multigraphs
such that $(A,f_{A})$ can be satisfied if and only if all these can be satisfied.

Therefore, the third stage of our algorithm will be to split our $2$-edge-connected components 
into $2$-vertex-connected blocks
(the decomposition is, in fact, unique),
and give each the appropriate discrepancy function.

\subsubsection*{STAGE 4}

Let $C$ be one of our $2$-vertex-connected blocks.
We wish to determine whether or not $(C,f_{C})$ can be satisfied. 
Analogously to Stages 1-3,
we shall split $C$ up into pieces at $2$-vertex-cuts.
However, unlike with these earlier stages,
this time if there exists a graph $M$ satisfying $(C,f_{C})$ 
there may be several different possibilities for how the edges of $M$ could interact with these pieces.
To keep track of this,
we shall define augmentations of the pieces
in such a way that $(C,f_{C})$ can be satisfied if and only if these augmentations can all be satisfied.

We will proceed iteratively.
At the start of each iteration,
we shall have a `blue' graph
(which will initially be $C$)
and an augmentation of it
(initially $(C,f_{C})$)
for which we want to determine satisfiability
(we shall refer to this augmentation as a `red' graph with a discrepancy function).
We will split our blue graph in two at a $2$-vertex-cut by breaking off a $3$-vertex-connected piece,
and we shall define augmentations of these two pieces
in such a way that the augmentation of the blue graph can be satisfied if and only if
the augmentations of the pieces can.
Lemma~\ref{lemma2} can then be used to determine satisfiability of the augmentation of the $3$-vertex-connected piece,
while the other piece and its augmentation can be used as the inputs for the next iteration.
The iterative loop terminates when the blue graph is itself $3$-vertex-connected.

We shall now give the full details:

\subsubsection*{Initialising}

Let us define our initial `blue graph', $B$, to be $C$,
let us also define our initial `red graph', $R$, to be $C$,
and let $R$ have discrepancy function $f_{R}=f_{C}$.
Note that $(R,f_{R})$ is an augmentation of $B$.
At the start of each iteration,
we will always have a blue planar multigraph with no cut-vertex,
and an augmentation of this consisting of a red graph and a discrepancy function.

\subsubsection*{The Iterative Loop}

Check if $B$ has any $2$-vertex-cuts.
If not, then we are done, 
since we can simply use Lemma~\ref{lemma2}.
Otherwise, let us find a minimal $2$-vertex-cut $\{ u,v \}$,
where we use `minimal' to mean that the component of smallest order in $B -u-v$ is minimal
over all possible $2$-vertex-cuts.

\begin{figure} [ht] 
\setlength{\unitlength}{0.9cm}
\begin{picture}(10,11.5)(1.375,-9.5)

\put(0.5,0.5){\oval(1,1)}
\put(3.5,0.5){\oval(1,1)}
\put(5.5,0.5){\oval(2,2)}
\put(0.5,0){\line(1,0){3}}
\put(0.5,1){\line(1,0){3}}
\put(2,0){\line(0,1){1}}
\put(5.5,1){\oval(7,1)[tl]}
\put(5.5,0){\oval(7,1)[bl]}
\put(2,0){\circle*{0.1}}
\put(2,1){\circle*{0.1}}
\put(1.8,1.15){$u$}
\put(1.8,-0.35){$v$}
\put(1.8,1.8){\large{$B$}}
\put(0.5,0){\circle*{0.1}}
\put(0.5,1){\circle*{0.1}}
\put(3.5,0){\circle*{0.1}}
\put(3.5,1){\circle*{0.1}}
\put(5.25,-0.5){\circle*{0.1}}
\put(5.25,1.5){\circle*{0.1}}

\put(0.5,-2.5){\oval(1,1)}
\put(3.5,-2.5){\oval(1,1)}
\put(5.5,-2.5){\oval(2,2)}
\put(0.3,-1.7){\large{$B_{1}$}}
\put(4,-1.5){\large{$B_{2}$}}
\put(0.5,-3){\circle*{0.1}}
\put(0.5,-2){\circle*{0.1}}
\put(3.5,-3){\circle*{0.1}}
\put(3.5,-2){\circle*{0.1}}
\put(5.25,-3.5){\circle*{0.1}}
\put(5.25,-1.5){\circle*{0.1}}

\put(9.25,0.5){\oval(1,1)}
\put(12.25,0.5){\oval(1,1)}
\put(14.25,0.5){\oval(2,2)}
\put(9.25,0){\line(1,0){3}}
\put(9.25,1){\line(1,0){3}}
\put(10.75,0){\line(0,1){1}}
\put(14.25,1){\oval(7,1)[tl]}
\put(14.25,0){\oval(7,1)[bl]}
\put(10.75,0){\circle*{0.1}}
\put(10.75,1){\circle*{0.1}}
\put(10.55,1.15){$u$}
\put(10.55,-0.35){$v$}
\put(10.6,1.8){\large{$R$}}
\put(10.1,1){\circle*{0.1}}
\put(10.1,-0.4){\line(0,1){0.8}}
\put(10.1,-0.4){\line(1,1){0.4}}
\put(10.1,-0.4){\line(-1,1){0.4}}
\put(10.1,0.4){\line(1,-1){0.4}}
\put(10.1,0.4){\line(-1,-1){0.4}}
\put(9.7,0){\circle*{0.1}}
\put(10.1,0){\circle*{0.1}}
\put(10.5,0){\circle*{0.1}}
\put(10.1,0.4){\circle*{0.1}}
\put(10.1,-0.4){\circle*{0.1}}
\put(10.75,0.5){\circle*{0.1}}
\put(12.25,1.5){\circle*{0.1}}
\put(12.25,-0.5){\circle*{0.1}}
\put(11.4,1){\circle*{0.1}}
\put(9.25,0){\circle*{0.1}}
\put(9.25,1){\circle*{0.1}}
\put(12.25,0){\circle*{0.1}}
\put(12.25,1){\circle*{0.1}}
\put(14,-0.5){\circle*{0.1}}
\put(14,1.5){\circle*{0.1}}

\put(9.25,-2.5){\oval(1,1)}
\put(12.25,-2.5){\oval(1,1)}
\put(14.25,-2.5){\oval(2,2)}
\put(9.35,-1.7){\large{$R_{1}$}}
\put(12.75,-1.3){\large{$R_{2}$}}
\put(10.1,-2){\circle*{0.1}}
\put(9.25,-2){\line(1,0){0.85}}
\put(9.25,-3){\line(1,0){1.25}}
\put(10.1,-3.4){\line(0,1){0.8}}
\put(10.1,-3.4){\line(1,1){0.4}}
\put(10.1,-3.4){\line(-1,1){0.4}}
\put(10.1,-2.6){\line(1,-1){0.4}}
\put(10.1,-2.6){\line(-1,-1){0.4}}
\put(9.7,-3){\circle*{0.1}}
\put(10.1,-3){\circle*{0.1}}
\put(10.5,-3){\circle*{0.1}}
\put(10.1,-2.6){\circle*{0.1}}
\put(10.1,-3.4){\circle*{0.1}}
\put(12.25,-1.5){\circle*{0.1}}
\put(12.25,-3.5){\circle*{0.1}}
\put(11.4,-2){\circle*{0.1}}
\put(11.4,-2){\line(1,0){0.85}}
\put(12.25,-1.5){\line(1,0){2}}
\put(12.25,-3.5){\line(1,0){2}}
\put(9.25,-3){\circle*{0.1}}
\put(9.25,-2){\circle*{0.1}}
\put(12.25,-3){\circle*{0.1}}
\put(12.25,-2){\circle*{0.1}}
\put(14,-3.5){\circle*{0.1}}
\put(14,-1.5){\circle*{0.1}}

\put(0.25,-5.5){\oval(1,1)}
\put(3.75,-5.5){\oval(1,1)}
\put(5.75,-5.5){\oval(2,2)}
\put(0.25,-6){\line(1,0){1.5}}
\put(0.25,-5){\line(1,0){1.5}}
\put(2.25,-6){\line(1,0){1.5}}
\put(2.25,-5){\line(1,0){1.5}}
\put(2.25,-6){\circle*{0.1}}
\put(2.25,-5){\circle*{0.1}}
\put(1.75,-6){\line(0,1){1}}
\put(5.75,-5){\oval(7,1)[tl]}
\put(5.75,-6){\oval(7,1)[bl]}
\put(1.75,-6){\circle*{0.1}}
\put(1.75,-5){\circle*{0.1}}
\put(2.05,-4.85){$u$}
\put(2.05,-6.35){$v$}
\put(1.65,-4.85){$u$}
\put(1.65,-6.35){$v$}
\put(4.25,-4.3){\large{$B_{2}^{\dag}$}}
\put(0.65,-4.5){\large{$B_{1}^{*}$}}
\put(0.25,-6){\circle*{0.1}}
\put(0.25,-5){\circle*{0.1}}
\put(3.75,-6){\circle*{0.1}}
\put(3.75,-5){\circle*{0.1}}
\put(5.5,-6.5){\circle*{0.1}}
\put(5.5,-4.5){\circle*{0.1}}

\put(9,-5.5){\oval(1,1)}
\put(12.5,-5.5){\oval(1,1)}
\put(14.5,-5.5){\oval(2,2)}
\put(9,-6){\line(1,0){1.5}}
\put(9,-5){\line(1,0){1.5}}
\put(11,-6){\line(1,0){1.5}}
\put(11,-5){\line(1,0){1.5}}
\put(10.5,-6){\line(0,1){1}}
\put(14.5,-5){\oval(7,1)[tl]}
\put(14.5,-6){\oval(7,1)[bl]}
\put(10.5,-6){\circle*{0.1}}
\put(10.5,-5){\circle*{0.1}}
\put(11,-6){\circle*{0.1}}
\put(11,-5){\circle*{0.1}}
\put(10.4,-4.85){$u$}
\put(10.4,-6.35){$v$}
\put(10.8,-4.85){$u$}
\put(10.8,-6.35){$v$}
\put(13,-4.3){\large{$R_{2}^{\dag}$}}
\put(9.85,-5){\circle*{0.1}}
\put(9.85,-6.4){\line(0,1){0.8}}
\put(9.85,-6.4){\line(1,1){0.4}}
\put(9.85,-6.4){\line(-1,1){0.4}}
\put(9.85,-5.6){\line(1,-1){0.4}}
\put(9.85,-5.6){\line(-1,-1){0.4}}
\put(9.45,-6){\circle*{0.1}}
\put(9.85,-6){\circle*{0.1}}
\put(10.25,-6){\circle*{0.1}}
\put(9.85,-5.6){\circle*{0.1}}
\put(9.85,-6.4){\circle*{0.1}}
\put(10.5,-5.5){\circle*{0.1}}
\put(12.5,-4.5){\circle*{0.1}}
\put(12.5,-6.5){\circle*{0.1}}
\put(11.65,-5){\circle*{0.1}}
\put(9.4,-4.5){\large{$R_{1}^{*}$}}
\put(9,-6){\circle*{0.1}}
\put(9,-5){\circle*{0.1}}
\put(12.5,-6){\circle*{0.1}}
\put(12.5,-5){\circle*{0.1}}
\put(14.25,-6.5){\circle*{0.1}}
\put(14.25,-4.5){\circle*{0.1}}

\put(9,-8.5){\oval(1,1)}
\put(12.5,-8.5){\oval(1,1)}
\put(14.5,-8.5){\oval(2,2)}
\put(9,-9){\line(1,0){1.5}}
\put(9,-8){\line(1,0){1.5}}
\put(11,-9){\line(1,0){1.5}}
\put(11,-8){\line(1,0){1.5}}
\put(11,-9){\line(0,1){1}}
\put(14.5,-8){\oval(7,1)[tl]}
\put(14.5,-9){\oval(7,1)[bl]}
\put(10.5,-9){\circle*{0.1}}
\put(10.5,-8){\circle*{0.1}}
\put(11,-9){\circle*{0.1}}
\put(11,-8){\circle*{0.1}}
\put(10.4,-7.85){$u$}
\put(10.4,-9.35){$v$}
\put(10.8,-7.85){$u$}
\put(10.8,-9.35){$v$}
\put(13,-7.3){\large{$R_{2}^{*}$}}
\put(9.85,-8){\circle*{0.1}}
\put(9.85,-9.4){\line(0,1){0.8}}
\put(9.85,-9.4){\line(1,1){0.4}}
\put(9.85,-9.4){\line(-1,1){0.4}}
\put(9.85,-8.6){\line(1,-1){0.4}}
\put(9.85,-8.6){\line(-1,-1){0.4}}
\put(9.45,-9){\circle*{0.1}}
\put(9.85,-9){\circle*{0.1}}
\put(10.25,-9){\circle*{0.1}}
\put(9.85,-8.6){\circle*{0.1}}
\put(9.85,-9.4){\circle*{0.1}}
\put(11,-8.5){\circle*{0.1}}
\put(12.5,-7.5){\circle*{0.1}}
\put(12.5,-9.5){\circle*{0.1}}
\put(11.65,-8){\circle*{0.1}}
\put(9.4,-7.5){\large{$R_{1}^{\dag}$}}
\put(9,-9){\circle*{0.1}}
\put(9,-8){\circle*{0.1}}
\put(12.5,-9){\circle*{0.1}}
\put(12.5,-8){\circle*{0.1}}
\put(14.25,-9.5){\circle*{0.1}}
\put(14.25,-7.5){\circle*{0.1}}

\put(0.25,-8.5){\oval(1,1)}
\put(3.75,-8.5){\oval(1,1)}
\put(5.75,-8.5){\oval(2,2)}
\put(0.25,-9){\line(1,0){1.5}}
\put(0.25,-8){\line(1,0){1.5}}
\put(2.25,-9){\line(1,0){1.5}}
\put(2.25,-8){\line(1,0){1.5}}
\put(2.25,-9){\circle*{0.1}}
\put(2.25,-8){\circle*{0.1}}
\put(2.25,-9){\line(0,1){1}}
\put(5.75,-8){\oval(7,1)[tl]}
\put(5.75,-9){\oval(7,1)[bl]}
\put(1.75,-9){\circle*{0.1}}
\put(1.75,-8){\circle*{0.1}}
\put(2.05,-7.85){$u$}
\put(2.05,-9.35){$v$}
\put(1.65,-7.85){$u$}
\put(1.65,-9.35){$v$}
\put(4.25,-7.3){\large{$B_{2}^{*}$}}
\put(0.65,-7.5){\large{$B_{1}^{\dag}$}}
\put(0.25,-9){\circle*{0.1}}
\put(0.25,-8){\circle*{0.1}}
\put(3.75,-9){\circle*{0.1}}
\put(3.75,-8){\circle*{0.1}}
\put(5.5,-9.5){\circle*{0.1}}
\put(5.5,-7.5){\circle*{0.1}}

\end{picture}
\caption{$\textrm{The planar multigraphs defined in the iterative loop of Stage 4.}$}
\label{D7}
\end{figure}
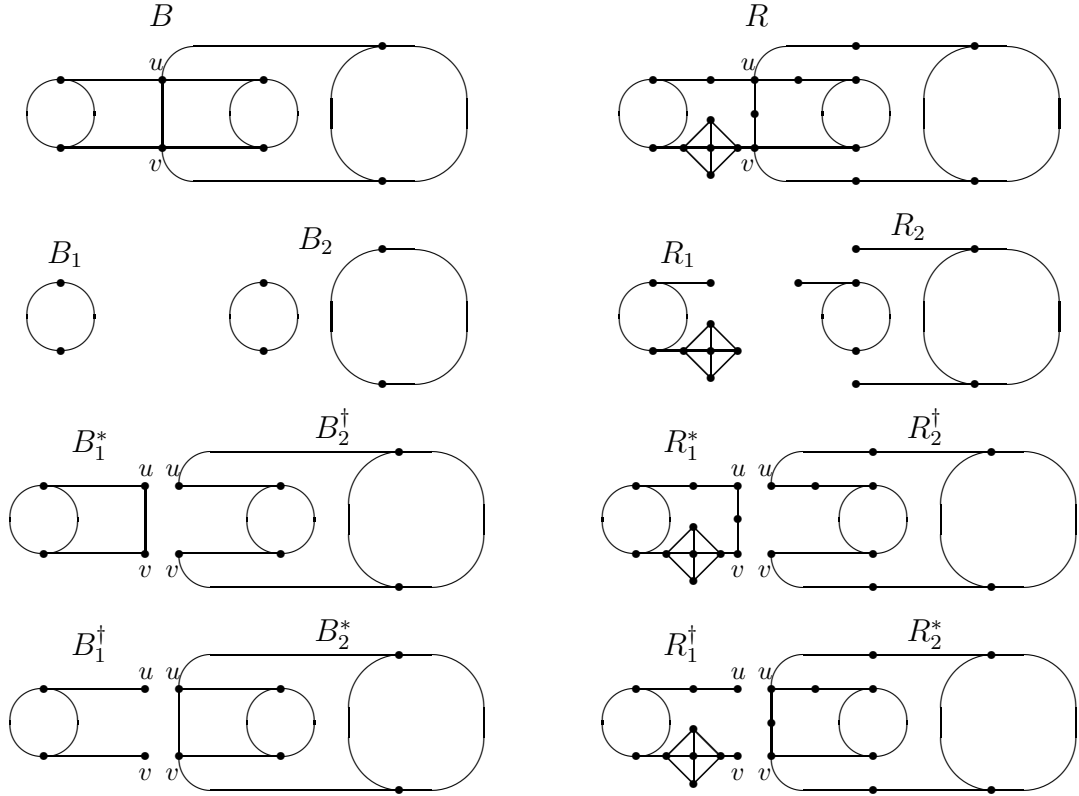

We shall now proceed to define several graphs based on the pieces of $B- u -v$
(these definitions are illustrated in Figure~\ref{D7}).
Let $B_{1}$ denote a component of smallest order in $B -u-v$,
let $B_{1}^{*}$ denote the graph induced by $V(B_{1}) \cup \{ u,v \}$ \label{b1star}
and let $B_{1}^{\dag}$ denote the graph obtained from $B_{1}^{*}$ by deleting any edges from $u$ to $v$.
Let $B_{2} = B \setminus B_{1}^{*}$,
let $B_{2}^{*} = B \setminus B_{1}$
and let $B_{2}^{\dag}$ denote the graph obtained from $B_{2}^{*}$ by deleting any edges from $u$ to $v$.
Let $R_{1}^{*},R_{2}^{*},R_{1}^{\dag}$ and $R_{2}^{\dag}$,
respectively,
denote the red versions of $B_{1}^{*},B_{2}^{*},B_{1}^{\dag}$ and~$B_{2}^{\dag}$
that follow `naturally' from $R$,
and let $R_{1} = R \setminus R_{2}^{*}$
and $R_{2} = R \setminus R_{1}^{*}$.

Let $u1$ denote the statement 
\begin{quote}
`$f_{R}(u)=0$ or there is only one edge in $B$ from $u$ to $B_{1}$' 
\end{quote}
(note that the latter clause of $u1$ implies $|B_{1}|=1$, by the minimality of $\{ u,v \}$,
but that it is not equivalent to this,
as we may have multi-edges).
It is important to note that the number of edges in $B$ from $u$ or $v$ to $B_{1}$
is exactly the same as the number of edges in $R$ from $u$ or $v$, respectively, to $R_{1}$
(and similarly for $B_{2}$ and $R_{2}$).
Thus, $u1$ is equivalent to the statement
`$f_{R}(u)=0$ or there is only one edge in $R$ from $u$ to $R_{1}$'.
Let $v1$ denote the analogous statement to $u1$ for $v$,
and let $u2$ and $v2$ denote the analogous statements for $B_{2}$.
Let $\overline{u1},\overline{v1},\overline{u2}$ and $\overline{v2}$
denote the complements of $u1,v1,u2$ and $v2$. 

Recall that we wish to split our graph in two at each iteration.
Note that if we have $\overline{u2}$, for example,
then $f_{R}(u) \geq 1$
and there are at least two edges in $R$ from $u$ to $R_{2}$,
so there may be several possibilities for where a graph satisfying $(R,f_{R})$
could have a new $u-R_{2}$ edge.
This could complicate matters,
causing an exponential blow-up in the running time,
unless we choose to split the graph in such a way that only the edges from $u$ to $R_{1}$
are important to the analysis.
Thus, our choice of how best to split the graph depends on which of the statements $u1,v1,u2$ and $v2$ are true,
and hence our next step is to divide our iterative loop into different cases based on this information.

\subsubsection*{Case (a): \boldmath{$u1 \land v1$}}

We shall first establish a couple of important facts,
before then splitting into two further subcases arising from parity issues.

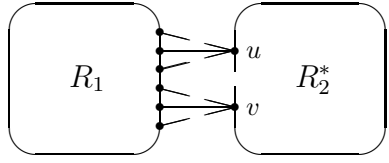
\begin{figure} [ht] 
\setlength{\unitlength}{1cm}
\begin{picture}(10,2)(-3.5,0)

\put(1,1.625){\oval(2,0.75)[t]}
\put(4,1.625){\oval(2,0.75)[t]}
\put(1,0.375){\oval(2,0.75)[b]}
\put(4,0.375){\oval(2,0.75)[b]}
\put(0,0.35){\line(0,1){1.3}}
\put(2,0.35){\line(0,1){1.3}}
\put(3,0.35){\line(0,1){0.55}}
\put(3,1.1){\line(0,1){0.55}}
\put(5,0.35){\line(0,1){1.3}}

\put(2,0.625){\line(1,0){1}}
\put(2,1.375){\line(1,0){1}}

\put(3,0.625){\circle*{0.1}}
\put(3,1.375){\circle*{0.1}}

\put(2,0.875){\circle*{0.1}}
\put(2,0.625){\circle*{0.1}}
\put(2,0.375){\circle*{0.1}}
\put(2,1.625){\circle*{0.1}}
\put(2,1.375){\circle*{0.1}}
\put(2,1.125){\circle*{0.1}}

\put(3.15,1.275){$u$}
\put(3.15,0.525){$v$}
\put(0.8,0.9){\large{$R_{1}$}}
\put(3.8,0.9){\large{$R_{2}^{*}$}}

\put(3,0.625){\line(-4,-1){0.4}}
\put(2.4,0.475){\line(-4,-1){0.4}}
\put(3,0.625){\line(-4,1){0.4}}
\put(2.4,0.775){\line(-4,1){0.4}}
\put(3,1.375){\line(-4,-1){0.4}}
\put(2.4,1.225){\line(-4,-1){0.4}}
\put(3,1.375){\line(-4,1){0.4}}
\put(2.4,1.525){\line(-4,1){0.4}}

\end{picture}
\caption{$\textrm{The planar multigraph $R$.}$}
\label{D8}
\end{figure}

Note that, by definition,
$B_{1}$ is connected.
Thus, $R_{1}$ must also be connected,
and so has to lie in a single face of $R_{2}^{*}$.
Hence, in any planar embedding $R$ must look as in Figure~\ref{D8},
where broken lines represent edges that may or may not exist
and where, without loss of generality,
we have drawn $R_{1}$ in the outside face of $R_{2}^{*}$.
Therefore, if a plane multigraph $M$ satisfies $(R,f_{R})$
then all edges in $E(M) \setminus~E(R)$ between $V(R_{1})$ and $V(R_{2}^{*})$
must lie within only two faces of the induced embedding of $R$
(since, by the definition of $u1$, $u$ can have more than one edge to $R_{1}$ only if $f(u)=0$,
and similarly for $v$).

Secondly,
since $f_{R}$ satisfies discrepancy parity,
note that $\sum_{x \in V(R_{1})} f_{R}(x)$ and $\sum_{x \in V \left( R_{2}^{*} \right) } f_{R}(x)$ 
must either both be odd or both be even.

\subsubsection*{Case (a)(i): 
\boldmath{$\sum_{x \in V(R_{1})} f_{R}(x)$} and \boldmath{$\sum_{x \in V \left( R_{2}^{*} \right) } f_{R}(x)$} both odd}

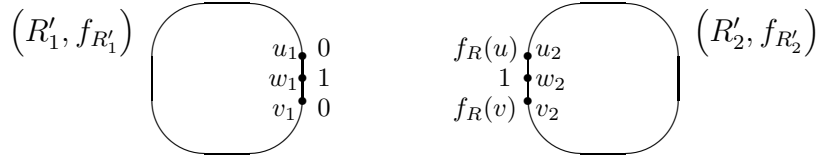
\begin{figure} [ht] 
\setlength{\unitlength}{1cm}
\begin{picture}(10,2)(-3.5,0)

\put(0,1){\oval(2,2)}
\put(5,1){\oval(2,2)}
\put(1,0.7){\circle*{0.1}}
\put(1,1){\circle*{0.1}}
\put(1,1.3){\circle*{0.1}}
\put(4,0.7){\circle*{0.1}}
\put(4,1){\circle*{0.1}}
\put(4,1.3){\circle*{0.1}}
\put(-2.9,1.5){\large{$\left( R_{1}^{\prime},f_{R_{1}^{\prime}} \right)$}}
\put(6.2,1.5){\large{$\left( R_{2}^{\prime},f_{R_{2}^{\prime}} \right)$}}
\put(0.6,1.3){$u_{1}$}
\put(0.55,0.9){$w_{1}$}
\put(0.6,0.5){$v_{1}$}
\put(4.1,1.3){$u_{2}$}
\put(4.1,0.9){$w_{2}$}
\put(4.1,0.5){$v_{2}$}
\put(1.2,1.3){$0$}
\put(1.2,0.9){$1$}
\put(1.2,0.5){$0$}
\put(3,1.3){$f_{R}(u)$}
\put(3.6,0.9){$1$}
\put(3,0.5){$f_{R}(v)$}

\end{picture}
\caption{$\textrm{The planar multigraphs $R_{1}^{\prime}$ and $R_{2}^{\prime}$, with their discrepancy functions.}$}
\label{D9}
\end{figure}

Let $B_{1}^{\prime} = B_{1}^{\dag} + uv$ and let $B_{2}^{\prime} = B_{2}^{*} + uv$
(so $uv$ will now be a multi-edge in $B_{2}^{\prime}$ if $uv \in E(B)$).
We shall now define an augmentation $(R_{1}^{\prime}, f_{R_{1}^{\prime}})$ of $B_{1}^{\prime}$
and an augmentation $(R_{2}^{\prime}, f_{R_{2}^{\prime}})$ of $B_{2}^{\prime}$
such that $(R,f_{R})$ can be satisfied if and only if
$\left( R_{1}^{\prime},f_{R_{1}^{\prime}} \right)$ and $\left( R_{2}^{\prime},f_{R_{2}^{\prime}} \right)$ 
can both be satisfied
(these new augmentations are illustrated in Figure~\ref{D9}).

Let $R_{1}^{\prime}$ be the graph formed from $R_{1}^{\dag}$
by relabelling $u$ and $v$ as $u_{1}$ and $v_{1}$, respectively,
and introducing a new vertex $w_{1}$ with edges to both $u_{1}$ and $v_{1}$.
Similarly, let $R_{2}^{\prime}$ be the graph formed from $R_{2}^{*}$
by relabelling $u$ and $v$ as $u_{2}$ and $v_{2}$, respectively,
and introducing a new vertex $w_{2}$ with edges to both $u_{2}$ and $v_{2}$.
Let $f_{R_{1}^{\prime}}$ be the discrepancy function on $R_{1}^{\prime}$ defined by setting
$f_{R_{1}^{\prime}}(u_{1})=f_{R_{1}^{\prime}}(v_{1})=0$,
$f_{R_{1}^{\prime}}(w_{1})=1$, and $f_{R_{1}^{\prime}}(x)=f_{R}(x)$~$\forall x \in V(R_{1})$.
Let $f_{R_{2}^{\prime}}$ be the discrepancy function on $R_{2}^{\prime}$ defined by setting
$f_{R_{2}^{\prime}}(u_{2})=f_{R}(u)$,
$f_{R_{2}^{\prime}}(v_{2})=f_{R}(v)$,
$f_{R_{2}^{\prime}}(w_{2})=1$, 
and $f_{R_{2}^{\prime}}(x)=f_{R}(x)$~$\forall x \in V(R_{2})$.
(Note that $f_{R_{1}^{\prime}}$ and $f_{R_{2}^{\prime}}$ are both valid discrepancy functions, 
since the discrepancy inequality is clearly satisfied by both
and discrepancy parity follows from the facts that
$\sum_{x \in V \left( R_{1}^{\prime} \right) } f_{R_{1}^{\prime}}(x) = \sum_{x \in V(R_{1})} f_{R}(x)+1$,
that
$\sum_{x \in V \left( R_{2}^{\prime} \right) } f_{R_{2}^{\prime}}(x) = \sum_{x \in V \left( R_{2}^{*} \right) } f_{R}(x)+1$
and that
$\sum_{x \in V(R_{1})} f_{R}(x)$ and $\sum_{x \in V \left( R_{2}^{*} \right) } f_{R}(x)$ are both odd).

\begin{Claim}
$(R,f_{R})$ can be satisfied if and only if
$\left( R_{1}^{\prime},f_{R_{1}^{\prime}} \right)$ and $\left( R_{2}^{\prime},f_{R_{2}^{\prime}} \right)$ 
can both be satisfied.
\end{Claim}
\textbf{Proof}
Suppose first that a plane multigraph $M$ satisfies $(R,f_{R})$.
Since $\sum_{x \in V(R_{1})} f_{R}(x)$ and $\sum_{x \in V \left( R_{2}^{*} \right) } f_{R}(x)$ are both odd,
there must be an odd number of edges in $E(M) \setminus E(R)$ between $V(R_{1})$ and $V(R_{2}^{*})$.
As already noted,
these edges must all lie within two faces of the embedding of $R$ induced from $M$.
Thus, one of these faces must have an odd number of new edges 
and the other must have an even number.
By pairing edges up,
as in the second half of Stage $3$,
we can hence obtain a planar multigraph satisfying $(R,f_{R})$ 
that has exactly one new edge between $V(R_{1})$ and $V(R_{2}^{*})$.
It is then easy to see that we can satisfy both 
$\left( R_{1}^{\prime},f_{R_{1}^{\prime}} \right)$ and $\left( R_{2}^{\prime},f_{R_{2}^{\prime}} \right)$.

Suppose next that 
$\left( R_{1}^{\prime},f_{R_{1}^{\prime}} \right)$ and $\left( R_{2}^{\prime},f_{R_{2}^{\prime}} \right)$ 
can both be satisfied,
by plane multigraphs $M_{R_{1}^{\prime}}$ and $M_{R_{2}^{\prime}}$ respectively,
and let the edges adjacent to $w_{1}$ in $E \left( M_{R_{1}^{\prime}} \right) \setminus E(R_{1}^{\prime})$ and 
$w_{2}$ in $E \left( M_{R_{2}^{\prime}} \right) \setminus E(R)$
be denoted by $e_{1}=z_{1}w_{1}$ and $e_{2}=z_{2}w_{2}$ respectively.
We may assume that $e_{1}$ is in the outside face of $M_{R_{1}^{\prime}}$.
Note that the edges $u_{1}w_{1}$ and $v_{1}w_{1}$ must then be in the outside face of $M_{R_{1}^{\prime}} - e_{1}$,
since these are the only edges incident to $w_{1}$ in $M_{R_{1}^{\prime}} - e_{1}$.
Hence,
by turning our drawing upside-down if necessary,
we may assume that $u_{1},w_{1}$ and $v_{1}$ 
are in clockwise order around this outer face of $M_{R_{1}^{\prime}} - e_{1}$,
and so $M_{R_{1}^{\prime}}$ is as shown in Figure~\ref{D10}
(where, without loss of generality, 
we have drawn $e_{1}$ so that $v_{1}$ is also in the outside face of $M_{R_{1}^{\prime}}$).
Similarly, we may assume that $M_{R_{2}^{\prime}}$ is also as shown in Figure~\ref{D10}.
It is then clear that we can delete $w_{1}$ and $w_{2}$,
`glue' $u_{1}$ to $u_{2}$ and $v_{1}$ to $v_{2}$
(i.e.~identify $u_{1}$ and $u_{2}$ and, separately, $v_{1}$ and $v_{2}$),
and insert the edge $z_{1}z_{2}$  
to obtain a plane multigraph $M_{R}$ that will satisfy $(R,f_{R})$
(note that it doesn't matter whether or not $z_{2} \in \{ u_{2},v_{2} \}$).
$\phantom{qwerty}
\setlength{\unitlength}{.25cm}
\begin{picture}(1,1)
\put(0,0){\line(1,0){1}}
\put(0,0){\line(0,1){1}}
\put(1,1){\line(-1,0){1}}
\put(1,1){\line(0,-1){1}}
\end{picture}$ \\

\begin{figure} [ht] 
\setlength{\unitlength}{0.75cm}
\begin{picture}(10,3)(-1.5,-0.5)

\put(1,1){\oval(2,2)}
\put(5,1){\oval(2,2)}
\put(2,0.7){\circle*{0.1}}
\put(2,1){\circle*{0.1}}
\put(2,1.3){\circle*{0.1}}
\put(4,0.7){\circle*{0.1}}
\put(4,1){\circle*{0.1}}
\put(4,1.3){\circle*{0.1}}
\put(-0.8,2.2){$M_{R_{1}^{\prime}}$}
\put(6.7,2.2){$M_{R_{2}^{\prime}}$}
\put(1.4,1.3){$u_{1}$}
\put(1.35,0.9){$w_{1}$}
\put(1.4,0.5){$v_{1}$}
\put(4.3,1.3){$u_{2}$}
\put(4.3,0.9){$w_{2}$}
\put(4.3,0.5){$v_{2}$}
\put(1,2){\circle*{0.1}}
\put(0.5,2.1){$z_{1}$}
\put(2,1.75){\oval(1,1.5)[r]}
\put(1.75,2){\oval(1.5,1)[tl]}
\put(4,1.75){\oval(1,1.5)[l]}
\put(5.75,0){\oval(1.5,1)[b]}
\put(6,1.75){\oval(1,1.5)[tr]}
\put(4,2.5){\line(1,0){2}}
\put(6.5,0){\line(0,1){2}}
\put(1.75,2.5){\line(1,0){0.25}}
\put(5,0){\circle*{0.1}}
\put(4.5,-0.3){$z_{2}$}
\put(1.65,2.6){$e_{1}$}
\put(4.9,2.6){$e_{2}$}

\put(7.5,1){\vector(1,0){1}}

\put(10.5,1){\oval(2,2)}
\put(11.5,0.7){\circle*{0.1}}
\put(11.5,1.3){\circle*{0.1}}
\put(10.5,2){\circle*{0.1}}
\put(11,1.2){$u$}
\put(11,0.6){$v$}
\put(10,2.1){$z_{1}$}
\put(11.25,2){\oval(1.5,1)[tl]}
\put(12.5,1){\oval(2,2)}
\put(13.25,0){\oval(1.5,1)[b]}
\put(13.5,1.75){\oval(1,1.5)[tr]}
\put(11.25,2.5){\line(1,0){2.25}}
\put(14,0){\line(0,1){2}}
\put(12.5,0){\circle*{0.1}}
\put(12,-0.3){$z_{2}$}
\put(8.9,2.2){$M_{R}$}

\end{picture}
\caption{$\textrm{Constructing a planar multigraph $M_{R}$ satisfying $(R,f_{R})$.}$}
\label{D10}
\end{figure}
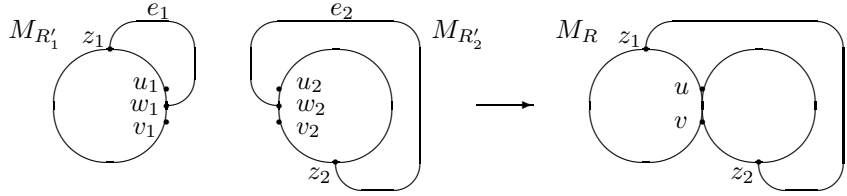

Recall that $B_{1}^{\prime} = B_{1}^{\dag} + uv$ 
and note that $B_{1}^{\prime}$ must not contain any $2$-vertex-cuts, by the minimality of $B_{1}$.
Thus, by Lemma~\ref{lemma2},
in $O \left( |B_{1}^{\prime}|^{2.5} \right)$ time
we can determine whether or not $\left( R_{1}^{\prime},f_{R_{1}^{\prime}} \right)$ can be satisfied.
If it cannot, we terminate the algorithm.
If it can, we return to the start of the iterative loop with $B_{2}^{\prime}$ as our new blue graph,
$R_{2}^{\prime}$ as our new red graph and $f_{R_{2}^{\prime}}$ as our new discrepancy function
(note that, as required,
$B_{2}^{\prime}$ does not contain a cut-vertex
since otherwise this would also be a cut-vertex in $B$
--- this property will be required for case~(b)).

\subsubsection*{Case (a)(ii): 
\boldmath{$\sum_{x \in V(R_{1})} f_{R}(x)$} and \boldmath{$\sum_{x \in V \left( R_{2}^{*} \right) } f_{R}(x)$} both even}

\begin{figure} [ht] 
\setlength{\unitlength}{1cm}
\begin{picture}(10,7)(-3,1)

\put(1,7.625){\oval(2,0.75)[t]}
\put(1,6.375){\oval(2,0.75)[b]}
\put(0,6.35){\line(0,1){1.3}}
\put(2,6.35){\line(0,1){1.3}}
\put(-1.9,7.5){\large{$\left( R_{1}^{\prime},f_{R_{1}^{\prime}} \right)$}}
\put(2,6.375){\circle*{0.1}}
\put(2,7.625){\circle*{0.1}}
\put(1.6,7.6){$u_{1}$}
\put(1.6,6.3){$v_{1}$}
\put(2.2,7.55){\footnotesize{$0$}}
\put(2.2,6.3){\footnotesize{$0$}}

\put(5,7.625){\oval(2,0.75)[t]}
\put(5,6.375){\oval(2,0.75)[b]}
\put(4,6.35){\line(0,1){1.3}}
\put(6,6.35){\line(0,1){1.3}}
\put(6.2,7.5){\large{$\left( R_{2}^{\prime},f_{R_{2}^{\prime}} \right)$}}
\put(4,6.375){\circle*{0.1}}
\put(4,7.625){\circle*{0.1}}
\put(4.2,7.6){$u_{2}$}
\put(4.2,6.3){$v_{2}$}
\put(3.1,7.6){\footnotesize{$f_{R}(u)$}}
\put(3.1,6.3){\footnotesize{$f_{R}(v)$}}

\put(1,5.125){\oval(2,0.75)[t]}
\put(1,3.875){\oval(2,0.75)[b]}
\put(0,3.85){\line(0,1){1.3}}
\put(2,3.85){\line(0,1){1.3}}
\put(-1.9,5){\large{$\left( R_{1}^{\prime\prime},f_{R_{1}^{\prime\prime}} \right)$}}
\put(2,3.875){\circle*{0.1}}
\put(2,5.125){\circle*{0.1}}
\put(1.6,5.1){$u_{1}$}
\put(1.6,3.8){$v_{1}$}
\put(1.6,4.5){\line(1,0){0.8}}
\put(1.6,4.5){\line(1,1){0.4}}
\put(1.6,4.5){\line(1,-1){0.4}}
\put(2.4,4.5){\line(-1,1){0.4}}
\put(2.4,4.5){\line(-1,-1){0.4}}
\put(2,4.5){\circle*{0.1}}
\put(1.6,4.5){\circle*{0.1}}
\put(2.4,4.5){\circle*{0.1}}
\put(2,4.1){\circle*{0.1}}
\put(2,4.9){\circle*{0.1}}
\put(1.3,4.4){\footnotesize{$1$}}
\put(2.5,4.4){\footnotesize{$1$}}
\put(2.1,4.5){\footnotesize{$0$}}
\put(2.1,3.7){\footnotesize{$0$}}
\put(2.1,5.1){\footnotesize{$0$}}
\put(2.1,4){\footnotesize{$0$}}
\put(2.1,4.85){\footnotesize{$0$}}

\put(5,5.125){\oval(2,0.75)[t]}
\put(5,3.875){\oval(2,0.75)[b]}
\put(4,3.85){\line(0,1){1.3}}
\put(6,3.85){\line(0,1){1.3}}
\put(6.2,5){\large{$\left( R_{2}^{\prime\prime},f_{R_{2}^{\prime\prime}} \right)$}}
\put(4,3.875){\circle*{0.1}}
\put(4,5.125){\circle*{0.1}}
\put(4.2,5.1){$u_{2}$}
\put(4.2,3.8){$v_{2}$}
\put(3.6,4.5){\line(1,0){0.8}}
\put(3.6,4.5){\line(1,1){0.4}}
\put(3.6,4.5){\line(1,-1){0.4}}
\put(4.4,4.5){\line(-1,1){0.4}}
\put(4.4,4.5){\line(-1,-1){0.4}}
\put(4,4.5){\circle*{0.1}}
\put(3.6,4.5){\circle*{0.1}}
\put(4.4,4.5){\circle*{0.1}}
\put(4,4.1){\circle*{0.1}}
\put(4,4.9){\circle*{0.1}}
\put(3.3,4.4){\footnotesize{$1$}}
\put(4.5,4.4){\footnotesize{$1$}}
\put(3.75,4.5){\footnotesize{$0$}}
\put(3.15,3.75){\footnotesize{$f_{R}(v)$}}
\put(3.15,5.15){\footnotesize{$f_{R}(u)$}}
\put(3.75,4){\footnotesize{$0$}}
\put(3.75,4.85){\footnotesize{$0$}}

\put(5,2.625){\oval(2,0.75)[t]}
\put(5,1.375){\oval(2,0.75)[b]}
\put(4,1.35){\line(0,1){1.3}}
\put(6,1.35){\line(0,1){1.3}}
\put(6.2,2.5){\large{$\left( R_{2}^{\prime\prime\prime},f_{R_{2}^{\prime\prime\prime}} \right)$}}
\put(4,1.375){\circle*{0.1}}
\put(4,2.625){\circle*{0.1}}
\put(4.2,2.6){$u_{2}$}
\put(4.2,1.3){$v_{2}$}
\put(4,2){\circle*{0.1}}
\put(4.2,1.9){$w$}
\put(3.1,2.6){\footnotesize{$f_{R}(u)$}}
\put(3.1,1.3){\footnotesize{$f_{R}(v)$}}
\put(3.7,1.9){\footnotesize{$2$}}

\end{picture}
\caption{$\textrm{The planar multigraphs 
$R_{1}^{\prime},R_{1}^{\prime\prime},R_{2}^{\prime},R_{2}^{\prime\prime}$ and $R_{2}^{\prime\prime\prime}$,
with their discrepancy}$} 
\textrm{functions.}
\label{D11}
\end{figure}
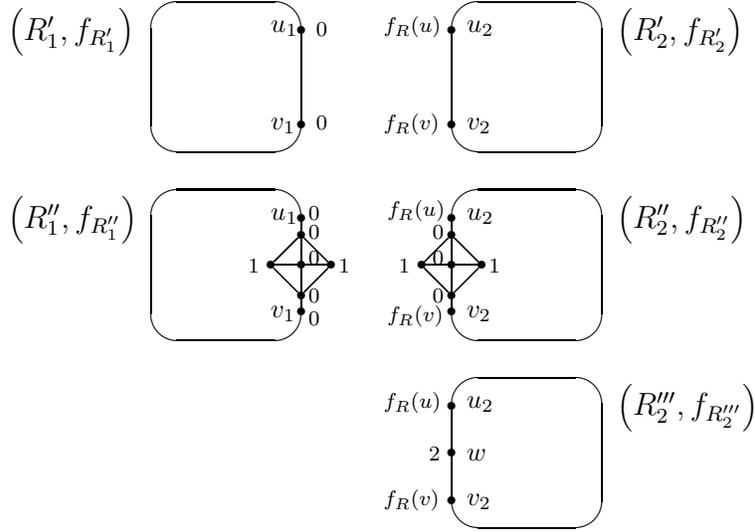

Again, we let $B_{1}^{\prime} = B_{1}^{\dag} + uv$ and $B_{2}^{\prime} = B_{2}^{*} + uv$.
This time,
we shall define augmentations 
$\left( R_{1}^{\prime},f_{R_{1}^{\prime}} \right)$ 
and $\left( R_{1}^{\prime\prime},f_{R_{1}^{\prime\prime}} \right)$ of $B_{1}^{\prime}$
and augmentations 
$\left( R_{2}^{\prime},f_{R_{2}^{\prime}} \right)$, $\left( R_{2}^{\prime\prime},f_{R_{2}^{\prime\prime}} \right)$
and $\left( R_{2}^{\prime\prime\prime},f_{R_{2}^{\prime\prime\prime}} \right)$ of $B_{2}^{\prime}$
(see Figure~\ref{D11})
such that
$(R,f_{R})$ can be satisfied if and only if: 
\begin{quote}
(1) $\left( R_{1}^{\prime},f_{R_{1}^{\prime}} \right)$ and $\left( R_{2}^{\prime},f_{R_{2}^{\prime}} \right)$ 
can both be satisfied, 
but $\left( R_{1}^{\prime\prime},f_{R_{1}^{\prime\prime}} \right)$~can't; \\ 
(2) $\left( R_{1}^{\prime\prime},f_{R_{1}^{\prime\prime}} \right)$ 
and $\left( R_{2}^{\prime\prime},f_{R_{2}^{\prime\prime}} \right)$ 
can both be satisfied, 
but $\left( R_{1}^{\prime},f_{R_{1}^{\prime}} \right)$~can't;~or \\
(3) $\left( R_{1}^{\prime},f_{R_{1}^{\prime}} \right),\left( R_{1}^{\prime\prime},f_{R_{1}^{\prime\prime}} \right)$ 
and $\left( R_{2}^{\prime\prime\prime},f_{R_{2}^{\prime\prime\prime}} \right)$ can all be satisfied. 
\end{quote}

Let $R_{1}^{\prime}$ be the graph formed from $R_{1}^{\dag}$ by relabelling $u$ and $v$ as $u_{1}$ and $v_{1}$,
respectively,
and inserting an edge between $u_{1}$ and $v_{1}$.
Let $f_{R_{1}^{\prime}}$ be the discrepancy function on $R_{1}^{\prime}$ defined by setting
$f_{R_{1}^{\prime}}(u_{1})=f_{R_{1}^{\prime}}(v_{1})=0$
and $f_{R_{1}^{\prime}}(x) =~f_{R}(x)$~$\forall x \in V(R_{1})$.
Let $R_{1}^{\prime\prime}$ be the graph formed from $R_{1}^{\prime}$ by placing a diamond on the $u_{1}v_{1}$ edge,
and let $f_{R_{1}^{\prime\prime}}$ be defined by setting
$f_{R_{1}^{\prime\prime}}(x)=f_{R_{1}^{\prime}}(x)$~$\forall x \in V(R_{1}^{\prime})$
and $f_{R_{1}^{\prime\prime}}(x) =~4 -\deg_{R_{1}^{\prime\prime}}(x)$~$\forall x \notin~V(R_{1}^{\prime})$.

Let $R_{2}^{\prime}$ be the graph formed from $R_{2}^{*}$ by relabelling $u$ and $v$ as $u_{2}$ and $v_{2}$,
respectively,
and inserting a new edge between $u_{1}$ and $v_{1}$
(so $u_{1}v_{1}$ will now be a multi-edge if $uv \in E(R)$).
Let $f_{R_{2}^{\prime}} = f_{R_{2}^{*}}$.
Let $R_{2}^{\prime\prime}$ be the graph formed from $R_{2}^{\prime}$ by
placing a diamond on the new $u_{2}v_{2}$ edge,
and let $f_{R_{2}^{\prime\prime}}$ be defined by setting
$f_{R_{2}^{\prime\prime}}(x)=f_{R_{2}^{\prime}}(x)$~$\forall x \in V(R_{2}^{\prime})$
and $f_{R_{2}^{\prime\prime}} = 4 - \deg_{R_{2}^{\prime\prime}}(x)$~$\forall x \notin V(R_{2}^{\prime})$.
Let $R_{2}^{\prime\prime\prime}$ be the graph formed from $R_{2}^{\prime}$
by instead subdividing the new $u_{2}v_{2}$ edge with a vertex $w$,
and let $f_{2}^{\prime\prime\prime}$ be defined by
$f_{2}^{\prime\prime\prime}(w)=2$
and $f_{2}^{\prime\prime\prime}(x) = f_{2}^{\prime}(x)$~$\forall x \in V(R_{2}^{\prime})$.

\begin{Claim}
 $(R,f_{R})$ can be satisfied if and only if
one of (1),(2) or (3) holds.
\end{Claim}
\textbf{Proof}
The `if' direction follows from a similar `gluing' argument as with case~(a)(i),
since we can again assume that the appropriate parts of our graphs are drawn in the outside face,
so we shall now proceed with proving the `only if' direction:

Suppose that a plane multigraph $M$ satisfies $(R,f_{R})$.
Since $\sum_{x \in V(R_{1})} f_{R}(x)$ and $\sum_{x \in V \left( R_{2}^{*} \right)} f_{R}(x)$ are both even,
there must be an even number of edges in $E(M) \setminus E(R)$ between $V(R_{1})$ and $V(R_{2}^{*})$.
As in case (a)(i),
these edges must all lie in two faces,
so we must either have an even number in both of these faces or an odd number in both.
By the same argument as with (a)(i),
we may in fact without loss of generality assume that 
there are either no new edges in both faces or exactly one in both.
In the former,
it is clear that we can satisfy both 
$\left( R_{1}^{\prime},f_{R_{1}^{\prime}} \right)$ and $\left( R_{2}^{\prime},f_{R_{2}^{\prime}} \right)$,
and in the latter it is clear that we can satisfy both 
$\left( R_{1}^{\prime\prime},f_{R_{1}^{\prime\prime}} \right)$ 
and $\left( R_{2}^{\prime\prime},f_{R_{2}^{\prime\prime}} \right)$.
Note that we can satisfy $\left( R_{2}^{\prime\prime\prime},f_{R_{2}^{\prime\prime\prime}} \right)$ if 
we can satisfy 
$\left( R_{2}^{\prime},f_{R_{2}^{\prime}} \right)$ or $\left( R_{2}^{\prime\prime},f_{R_{2}^{\prime\prime}} \right)$.
Thus, we can either satisfy 
$\left( R_{1}^{\prime},f_{R_{1}^{\prime}} \right)$,$\left( R_{2}^{\prime},f_{R_{2}^{\prime}} \right)$
and $\left( R_{2}^{\prime\prime\prime},f_{R_{2}^{\prime\prime\prime}} \right)$,
or $\left( R_{1}^{\prime\prime},f_{R_{1}^{\prime\prime}} \right)$,
$\left( R_{2}^{\prime\prime},f_{R_{2}^{\prime\prime}} \right)$ 
and $\left( R_{2}^{\prime\prime\prime},f_{R_{2}^{\prime\prime\prime}} \right)$.
In the first case, either (1) or (3) must hold,
and in the second case either (2) or (3) must hold.
$\phantom{qwerty}
\setlength{\unitlength}{.25cm}
\begin{picture}(1,1)
\put(0,0){\line(1,0){1}}
\put(0,0){\line(0,1){1}}
\put(1,1){\line(-1,0){1}}
\put(1,1){\line(0,-1){1}}
\end{picture}$ \\

We have now shown that $(R,f_{R})$ can be satisfied if and only if (1),(2) or (3) hold.
As in case (a)(i),
we can use Lemma~\ref{lemma2} to determine in $O \left( |B_{1}^{\prime}|^{2.5} \right)$ time whether 
$\left( R_{1}^{\prime},f_{R_{1}^{\prime}} \right)$ and $\left( R_{1}^{\prime\prime},f_{R_{1}^{\prime\prime}} \right)$ 
can be satisfied.
If neither can be satisfied,
we terminate the algorithm.
If at least one can be satisfied,
then we return to the start of the iterative loop with $B_{2}^{\prime}$ as our new blue graph 
and either 
$(R_{2}^{\prime\prime\prime},f_{R_{2}^{\prime\prime\prime}})$,
$(R_{2}^{\prime},f_{R_{2}^{\prime}})$ or $(R_{2}^{\prime\prime},f_{R_{2}^{\prime\prime}})$
as our augmentation,
according to whether both
$\left( R_{1}^{\prime},f_{R_{1}^{\prime}} \right)$ and $\left( R_{1}^{\prime\prime},f_{R_{1}^{\prime\prime}} \right)$,
just $\left( R_{1}^{\prime},f_{R_{1}^{\prime}} \right)$,
or just $\left( R_{1}^{\prime\prime},f_{R_{1}^{\prime\prime}} \right)$
can be satisfied, respectively.

\subsubsection*{Case (b): \boldmath{$(\overline{u1} \lor \overline{v1}) \land u2 \land v2$}}

We shall again start with some groundwork on the structure of $R$,
analogously to case (a),
before splitting into subcases.

Since $\overline{u1} \lor \overline{v1}$ holds,
we can't have $f(u) = f(v) =0$.
Thus, since $u2 \land v2$ also holds,
it must be that either $u$ or $v$ has only one edge to $B_{2}$.
Hence, since $B$ contains no cut-vertices,
it must be that $B_{2}$ is connected.
Therefore, $R_{2}$ must also be connected 
and so must lie in a single face of $R_{1}^{*}$.
Hence, we may proceed in a similar way to case (a),
but this time we will split into subcases depending on the parity of $R_{1}^{*}$ and $R_{2}$,
rather than $R_{1}$ and $R_{2}^{*}$.

\subsubsection*{Case (b)(i): 
\boldmath{$\sum_{x \in V \left( R_{1}^{*} \right) } f_{R}(x)$} and \boldmath{$\sum_{x \in V(R_{2})} f_{R}(x)$} both odd}

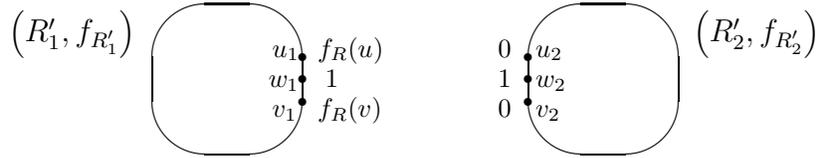
\begin{figure} [ht] 
\setlength{\unitlength}{1cm}
\begin{picture}(10,2)(-3.5,0)

\put(0,1){\oval(2,2)}
\put(5,1){\oval(2,2)}
\put(1,0.7){\circle*{0.1}}
\put(1,1){\circle*{0.1}}
\put(1,1.3){\circle*{0.1}}
\put(4,0.7){\circle*{0.1}}
\put(4,1){\circle*{0.1}}
\put(4,1.3){\circle*{0.1}}
\put(-2.9,1.5){\large{$\left( R_{1}^{\prime},f_{R_{1}^{\prime}} \right)$}}
\put(6.2,1.5){\large{$\left( R_{2}^{\prime},f_{R_{2}^{\prime}} \right)$}}
\put(0.6,1.3){$u_{1}$}
\put(0.55,0.9){$w_{1}$}
\put(0.6,0.5){$v_{1}$}
\put(4.1,1.3){$u_{2}$}
\put(4.1,0.9){$w_{2}$}
\put(4.1,0.5){$v_{2}$}
\put(1.2,1.3){$f_{R}(u)$}
\put(1.3,0.9){$1$}
\put(1.2,0.5){$f_{R}(v)$}
\put(3.6,1.3){$0$}
\put(3.6,0.9){$1$}
\put(3.6,0.5){$0$}

\end{picture}
\caption{$\textrm{The planar multigraphs $R_{1}^{\prime}$ and $R_{2}^{\prime}$, with their discrepancy functions.}$}
\label{D12}
\end{figure}

This time,
we let $B_{1}^{\prime} = B_{1}^{*} + uv$
(so $uv$ will be a multi-edge in $B_{1}^{\prime}$ if $uv \in E(B)$)
and let $B_{2}^{\prime} = B_{2}^{\dag} + uv$.
We will define augmentations
$(R_{1}^{\prime}, f_{R_{1}^{\prime}})$ of $B_{1}^{\prime}$
and $(R_{2}^{\prime}, f_{R_{2}^{\prime}})$ of $B_{2}^{\prime}$
(see Figure~\ref{D12})
such that $(R,f_{R})$ can be satisfied if and only if
$\left( R_{1}^{\prime},f_{R_{1}^{\prime}} \right)$ and $\left( R_{2}^{\prime},f_{R_{2}^{\prime}} \right)$ 
can both be satisfied.

Let $R_{1}^{\prime}$ be the graph formed from $R_{1}^{*}$
by relabelling $u$ and $v$ as $u_{1}$ and $v_{1}$, respectively,
and introducing a new vertex $w_{1}$ with edges to both $u_{1}$ and $v_{1}$.
Similarly, let $R_{2}^{\prime}$ be the graph formed from $R_{2}^{\dag}$
by relabelling $u$ and $v$ as $u_{2}$ and~$v_{2}$, respectively,
and introducing a new vertex $w_{2}$ with edges to both $u_{2}$ and $v_{2}$.
Let $f_{R_{1}^{\prime}}$ be the discrepancy function on $R_{1}^{\prime}$ defined by setting
$f_{R_{1}^{\prime}}(u_{1})=f_{R}(u)$,
$f_{R_{1}^{\prime}}(v_{1})=f_{R}(v)$,
$f_{R_{1}^{\prime}}(w_{1})=1$, 
and $f_{R_{1}^{\prime}}(x)=f_{R}(x)$~$\forall x \in V(R_{1})$.
Let $f_{R_{2}^{\prime}}$ be the discrepancy function on $R_{2}^{\prime}$ defined by setting
$f_{R_{2}^{\prime}}(u_{2})=f_{R_{2}^{\prime}}(v_{2})=0$,
$f_{R_{2}^{\prime}}(w_{2})=1$, and $f_{R_{2}^{\prime}}(x)=f_{R}(x)$~$\forall x \in V(R_{2})$.

The proof that $(R,f_{R})$ may be satisfied if and only if 
both $\left( R_{1}^{\prime},f_{R_{1}^{\prime}} \right)$ and $\left( R_{2}^{\prime},f_{R_{2}^{\prime}} \right)$ 
may be satisfied is as with case (a)(i).
Again,
we can determine in $O \left( |B_{1}^{\prime}|^{2.5} \right)$ time 
whether or not $\left( R_{1}^{\prime},f_{R_{1}^{\prime}} \right)$ can be satisfied,
and if so we return to the start of the iterative loop with $B_{2}^{\prime}$ as our new blue graph,
$R_{2}^{\prime}$ as our new red graph 
and $f_{R_{2}^{\prime}}$ as our new discrepancy function. 
Otherwise, we terminate the algorithm.

\subsubsection*{Case (b)(ii): 
\boldmath{$\sum_{x \in V \left( R_{1}^{*} \right)} f_{R}(x)$} and \boldmath{$\sum_{x \in V(R_{2})} f_{R}(x)$} both even}

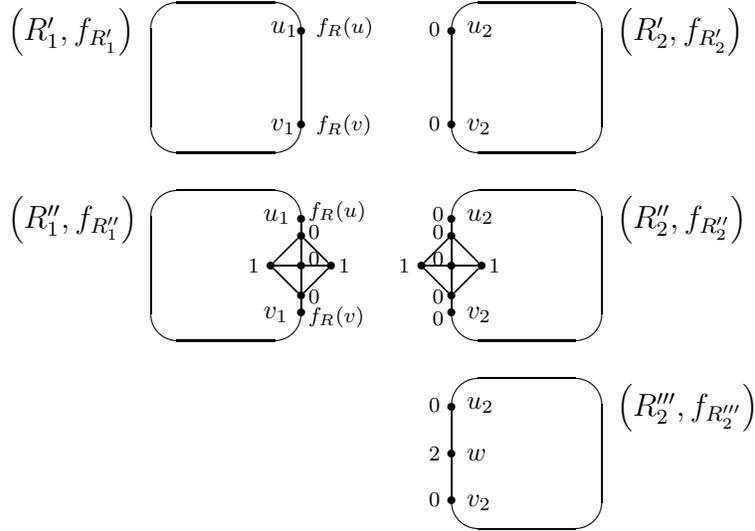
\begin{figure} [ht] 
\setlength{\unitlength}{1cm}
\begin{picture}(10,7)(-3,1)

\put(1,7.625){\oval(2,0.75)[t]}
\put(1,6.375){\oval(2,0.75)[b]}
\put(0,6.35){\line(0,1){1.3}}
\put(2,6.35){\line(0,1){1.3}}
\put(-1.9,7.5){\large{$\left( R_{1}^{\prime},f_{R_{1}^{\prime}} \right)$}}
\put(2,6.375){\circle*{0.1}}
\put(2,7.625){\circle*{0.1}}
\put(1.6,7.6){$u_{1}$}
\put(1.6,6.3){$v_{1}$}
\put(2.2,7.6){\footnotesize{$f_{R}(u)$}}
\put(2.2,6.3){\footnotesize{$f_{R}(v)$}}

\put(5,7.625){\oval(2,0.75)[t]}
\put(5,6.375){\oval(2,0.75)[b]}
\put(4,6.35){\line(0,1){1.3}}
\put(6,6.35){\line(0,1){1.3}}
\put(6.2,7.5){\large{$\left( R_{2}^{\prime},f_{R_{2}^{\prime}} \right)$}}
\put(4,6.375){\circle*{0.1}}
\put(4,7.625){\circle*{0.1}}
\put(4.2,7.6){$u_{2}$}
\put(4.2,6.3){$v_{2}$}
\put(3.7,7.55){\footnotesize{$0$}}
\put(3.7,6.3){\footnotesize{$0$}}

\put(1,5.125){\oval(2,0.75)[t]}
\put(1,3.875){\oval(2,0.75)[b]}
\put(0,3.85){\line(0,1){1.3}}
\put(2,3.85){\line(0,1){1.3}}
\put(-1.9,5){\large{$\left( R_{1}^{\prime\prime},f_{R_{1}^{\prime\prime}} \right)$}}
\put(2,3.875){\circle*{0.1}}
\put(2,5.125){\circle*{0.1}}
\put(1.5,5.1){$u_{1}$}
\put(1.5,3.8){$v_{1}$}
\put(1.6,4.5){\line(1,0){0.8}}
\put(1.6,4.5){\line(1,1){0.4}}
\put(1.6,4.5){\line(1,-1){0.4}}
\put(2.4,4.5){\line(-1,1){0.4}}
\put(2.4,4.5){\line(-1,-1){0.4}}
\put(2,4.5){\circle*{0.1}}
\put(1.6,4.5){\circle*{0.1}}
\put(2.4,4.5){\circle*{0.1}}
\put(2,4.1){\circle*{0.1}}
\put(2,4.9){\circle*{0.1}}
\put(1.3,4.4){\footnotesize{$1$}}
\put(2.5,4.4){\footnotesize{$1$}}
\put(2.1,4.5){\footnotesize{$0$}}
\put(2.1,3.75){\footnotesize{$f_{R}(v)$}}
\put(2.1,5.15){\footnotesize{$f_{R}(u)$}}
\put(2.1,4){\footnotesize{$0$}}
\put(2.1,4.85){\footnotesize{$0$}}

\put(5,5.125){\oval(2,0.75)[t]}
\put(5,3.875){\oval(2,0.75)[b]}
\put(4,3.85){\line(0,1){1.3}}
\put(6,3.85){\line(0,1){1.3}}
\put(6.2,5){\large{$\left( R_{2}^{\prime\prime},f_{R_{2}^{\prime\prime}} \right)$}}
\put(4,3.875){\circle*{0.1}}
\put(4,5.125){\circle*{0.1}}
\put(4.2,5.1){$u_{2}$}
\put(4.2,3.8){$v_{2}$}
\put(3.6,4.5){\line(1,0){0.8}}
\put(3.6,4.5){\line(1,1){0.4}}
\put(3.6,4.5){\line(1,-1){0.4}}
\put(4.4,4.5){\line(-1,1){0.4}}
\put(4.4,4.5){\line(-1,-1){0.4}}
\put(4,4.5){\circle*{0.1}}
\put(3.6,4.5){\circle*{0.1}}
\put(4.4,4.5){\circle*{0.1}}
\put(4,4.1){\circle*{0.1}}
\put(4,4.9){\circle*{0.1}}
\put(3.3,4.4){\footnotesize{$1$}}
\put(4.5,4.4){\footnotesize{$1$}}
\put(3.75,4.5){\footnotesize{$0$}}
\put(3.75,3.7){\footnotesize{$0$}}
\put(3.75,5.1){\footnotesize{$0$}}
\put(3.75,4){\footnotesize{$0$}}
\put(3.75,4.85){\footnotesize{$0$}}

\put(5,2.625){\oval(2,0.75)[t]}
\put(5,1.375){\oval(2,0.75)[b]}
\put(4,1.35){\line(0,1){1.3}}
\put(6,1.35){\line(0,1){1.3}}
\put(6.2,2.5){\large{$\left( R_{2}^{\prime\prime\prime},f_{R_{2}^{\prime\prime\prime}} \right)$}}
\put(4,1.375){\circle*{0.1}}
\put(4,2.625){\circle*{0.1}}
\put(4.2,2.6){$u_{2}$}
\put(4.2,1.3){$v_{2}$}
\put(4,2){\circle*{0.1}}
\put(4.2,1.9){$w$}
\put(3.7,2.55){\footnotesize{$0$}}
\put(3.7,1.3){\footnotesize{$0$}}
\put(3.7,1.9){\footnotesize{$2$}}

\end{picture}
\caption{$\textrm{The planar multigraphs $R_{1}^{\prime},R_{1}^{\prime\prime},R_{2}^{\prime},R_{2}^{\prime\prime}$
and $R_{2}^{\prime\prime\prime}$, with their discrepancy}$}
\textrm{functions.}
\label{D13}
\end{figure}

Again, we let $B_{1}^{\prime} = B_{1}^{*} + uv$
and $B_{2}^{\prime} = B_{2}^{\dag} + uv$.
This time, as with case (a)(ii),
we shall define augmentations 
$\left( R_{1}^{\prime},f_{R_{1}^{\prime}} \right)$ 
and $\left( R_{1}^{\prime\prime},f_{R_{1}^{\prime\prime}} \right)$ of $B_{1}^{\prime}$
and augmentations 
$\left( R_{2}^{\prime},f_{R_{2}^{\prime}} \right), \left( R_{2}^{\prime\prime},f_{R_{2}^{\prime\prime}} \right)$
and $\left( R_{2}^{\prime\prime\prime},f_{R_{2}^{\prime\prime\prime}} \right)$ of $B_{2}^{\prime}$
(see Figure~\ref{D13})
such that
$(R,f_{R})$ can be satisfied if and only if: 
\begin{quote}
(1) $\left( R_{1}^{\prime},f_{R_{1}^{\prime}} \right)$ and $\left( R_{2}^{\prime},f_{R_{2}^{\prime}} \right)$ 
can both be satisfied, 
but $\left( R_{1}^{\prime\prime},f_{R_{1}^{\prime\prime}} \right)$~can't; \\ 
(2) $\left( R_{1}^{\prime\prime},f_{R_{1}^{\prime\prime}} \right)$ 
and $\left( R_{2}^{\prime\prime},f_{R_{2}^{\prime\prime}} \right)$ 
can both be satisfied, 
but $\left( R_{1}^{\prime},f_{R_{1}^{\prime}} \right)$~can't;~or \\
(3) $\left( R_{1}^{\prime},f_{R_{1}^{\prime}} \right),\left( R_{1}^{\prime\prime},f_{R_{1}^{\prime\prime}} \right)$ 
and $\left( R_{2}^{\prime\prime\prime},f_{R_{2}^{\prime\prime\prime}} \right)$ can all be satisfied.
\end{quote}

Let $R_{1}^{\prime}$ be the graph formed from $R_{1}^{*}$ by relabelling $u$ and $v$ as $u_{1}$ and $v_{1}$,
respectively,
and inserting an edge between $u_{1}$ and $v_{1}$
(so $u_{1}v_{1}$ will now be a multi-edge if $uv \in E(R)$).
Let $f_{R_{1}^{\prime}}$ be defined by setting $f_{R_{1}^{\prime}}(x) = f_{R}(x)$~$\forall x \in~V(R_{1}^{\prime})$.
Let $R_{1}^{\prime\prime}$ be the graph formed from $R_{1}^{\prime}$ by placing a diamond on the $u_{1}v_{1}$ edge,
and let $f_{R_{1}^{\prime\prime}}$ be defined by setting
$f_{R_{1}^{\prime\prime}}(x)=f_{R_{1}^{\prime}}(x)$~$\forall x \in V(R_{1}^{\prime})$
and $f_{R_{1}^{\prime\prime}}(x) = 4 -\deg_{R_{1}^{\prime\prime}}(x)$~$\forall x \notin V(R_{1}^{\prime})$.

Let $R_{2}^{\prime}$ be the graph formed from $R_{2}^{\dag}$ by relabelling $u$ and $v$ as $u_{2}$ and~$v_{2}$,
respectively,
and inserting a new edge between $u_{2}$ and $v_{2}$.
Let $f_{R_{2}^{\prime}}$ be the discrepancy function on $R_{2}^{\prime}$ defined by setting
$f_{R_{2}^{\prime}}(u_{2})=f_{R_{2}^{\prime}}(v_{2})=0$
and $f_{R_{2}^{\prime}}(x) =~f_{R}(x)$~$\forall x \in V(R_{2})$.
Let $R_{2}^{\prime\prime}$ be the graph formed from $R_{2}^{\prime}$ by
placing a diamond on the new $u_{2}v_{2}$ edge,
and let $f_{R_{2}^{\prime\prime}}$ be defined by setting
$f_{R_{2}^{\prime\prime}}(x)=~f_{R_{2}^{\prime}}(x)$~$\forall x \in V(R_{2}^{\prime})$
and $f_{R_{2}^{\prime\prime}}(x) = 4 - \deg_{R_{2}^{\prime\prime}}(x)$~$\forall x \notin V(R_{2}^{\prime})$.
Let $R_{2}^{\prime\prime\prime}$ be the graph formed from $R_{2}^{\prime}$
by instead subdividing the new $u_{2}v_{2}$ edge with a vertex $w$,
and let $f_{R_{2}}^{\prime\prime\prime}$ be defined by
$f_{R_{2}}^{\prime\prime\prime}(w)=2$
and $f_{R_{2}}^{\prime\prime\prime}(x) =~f_{R_{2}}^{\prime}(x)$~$\forall x \in~V(R_{2}^{\prime})$.

The proof that $(R,f_{R})$ can be satisfied if and only if (1),(2) or (3) hold is as with case (a)(ii).
Again, we can determine in $O \left( |B_{1}^{\prime}|^{2.5} \right)$ time whether or not 
$\left( R_{1}^{\prime},f_{R_{1}^{\prime}} \right)$ and $\left( R_{1}^{\prime\prime},f_{R_{1}^{\prime\prime}} \right)$ 
can be satisfied,
and if at least one can
then we return to the start of the iterative loop with $B_{2}^{\prime}$ as our new blue graph 
and either 
$(R_{2}^{\prime\prime\prime},f_{R_{2}^{\prime\prime\prime}})$,
$(R_{2}^{\prime},f_{R_{2}^{\prime}})$ or $(R_{2}^{\prime\prime},f_{R_{2}^{\prime\prime}})$ as our augmentation,
according to whether both 
$\left( R_{1}^{\prime},f_{R_{1}^{\prime}} \right)$ and $\left( R_{1}^{\prime\prime},f_{R_{1}^{\prime\prime}} \right)$,
just $\left( R_{1}^{\prime},f_{R_{1}^{\prime}} \right)$,
or just $\left( R_{1}^{\prime\prime},f_{R_{1}^{\prime\prime}} \right)$
can be satisfied, respectively. 
If neither $\left( R_{1}^{\prime},f_{R_{1}^{\prime}} \right)$
nor $\left( R_{1}^{\prime\prime},f_{R_{1}^{\prime\prime}} \right)$ can be satisfied,
we terminate the algorithm.

\subsubsection*{Case (c): \boldmath{$(\overline{u1} \lor \overline{v1}) \land (\overline{u2} \lor \overline{v2})$}}

We will now deal with the remaining case,
which will follow from a detailed investigated of the properties that are forced upon us 
if $(\overline{u1} \lor \overline{v1}) \land (\overline{u2} \lor \overline{v2})$ holds.

Recall that if we have $\overline{u1}$,
then by definition $f_{R}(u) \geq 1$ and $u$ has at least two edges to $R_{1}$,
so $u$ must have only one edge to $R_{2}$,
and hence we have $u2$.
Similarly,
$\overline{v1} \Rightarrow v2$, $\overline{u2} \Rightarrow u1$ and $\overline{v2} \Rightarrow v1$.
Thus, the only possibilities are $u1 \land \overline{u2} \land \overline{v1} \land v2$
and $\overline{u1} \land u2 \land v1 \land \overline{v2}$.
By swapping $u$ and $v$ if necessary,
we can without loss of generality assume that we have the former.

\begin{figure} [ht] 
\setlength{\unitlength}{1cm}
\begin{picture}(10,2)(-3,2.5)

\put(4,4.125){\oval(2,0.75)[t]}
\put(4,2.875){\oval(2,0.75)[b]}
\put(3,2.85){\line(0,1){1.3}}
\put(5,2.85){\line(0,1){1.3}}
\put(2,3.9){\line(1,0){1}}
\put(2,3.9){\line(5,1){1}}
\put(2,3.1){\line(1,0){1}}
\put(1,3.5){\circle*{0.1}}
\put(2,3.1){\circle*{0.1}}
\put(2,3.9){\circle*{0.1}}
\put(3,3.1){\circle*{0.1}}
\put(1,3.5){\line(5,2){1}}
\put(2,3.5){\oval(2,0.8)[bl]}
\put(1,3.1){\oval(2,0.8)[tr]}
\put(3.8,3.4){\large{$B_{2}$}}
\put(0.8,3.7){\large{$B_{1}$}}
\put(1.9,4){$u$}
\put(1.9,2.8){$v$}
\put(2.8,2.8){$x$}
\put(3,3.9){\circle*{0.1}}
\put(3,4.1){\circle*{0.1}}

\end{picture}
\caption{$\textrm{The structure of $B$ in case (c).}$}
\label{D14a}
\end{figure}
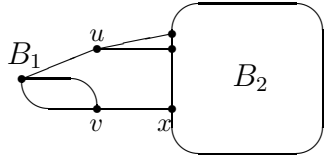

Note that the only way to obtain $u1 \land \overline{u2}$
is to have exactly one edge in $B$ from $u$ to $B_{1}$
(or, equivalently, exactly one edge in $R$ from $u$ to $R_{1}$),
exactly two edges in $B$ from $u$ to $B_{2}$,
no edges in $B$ from $u$ to $v$,
and $f_{R}(u)=1$.
Similarly, we must have exactly one edge in $B$ from $v$ to $B_{2}$,
exactly two edges in $B$ from $v$ to $B_{1}$,
and $f_{R}(v)=1$.
Note also that we must have $|B_{1}|=1$,
since otherwise the minimality of $B_{1}$ would imply that $u$ and $v$ 
would both have to have at least two edges in $B$ to $B_{1}$,
which would in turn imply that we would have to have $u2 \land v2$.
Thus, $B$ must be as shown in Figure~\ref{D14a}.

If $|B_{2}|=1$,
then $|R|$ is bounded by a constant
and so we can determine 
the satisfiability of $(R,f_{R})$ in $O(1)$ time
(simply by checking all graphs with $|R|$ vertices to see if any of these do satisfy $(R,f_{R})$).

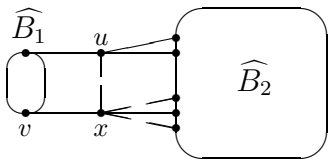
\begin{figure} [ht] 
\setlength{\unitlength}{1cm}
\begin{picture}(10,2)(-3,0)

\put(4,1.625){\oval(2,0.75)[t]}
\put(4,0.375){\oval(2,0.75)[b]}
\put(3,0.35){\line(0,1){1.3}}
\put(5,0.35){\line(0,1){1.3}}
\put(2,1.4){\line(1,0){1}}
\put(2,1.4){\line(5,1){1}}
\put(2,0.6){\line(1,0){1}}
\put(2,0.6){\line(5,1){0.4}}
\put(2.6,0.72){\line(5,1){0.4}}
\put(2,0.6){\line(5,-1){0.4}}
\put(2.6,0.48){\line(5,-1){0.4}}
\put(2,0.6){\line(0,1){0.36}}
\put(2,1.08){\line(0,1){0.36}}
\put(1,1.4){\line(1,0){1}}
\put(1,0.6){\line(1,0){1}}
\put(1,1){\oval(0.5,0.8)}
\put(1,0.6){\circle*{0.1}}
\put(1,1.4){\circle*{0.1}}
\put(2,0.6){\circle*{0.1}}
\put(2,1.4){\circle*{0.1}}
\put(3.8,0.9){\large{$\widehat{B_{2}}$}}
\put(0.8,1.6){\large{$\widehat{B_{1}}$}}
\put(1.9,1.5){$u$}
\put(1.9,0.3){$x$}
\put(0.9,0.3){$v$}
\put(3,0.4){\circle*{0.1}}
\put(3,0.6){\circle*{0.1}}
\put(3,0.8){\circle*{0.1}}
\put(3,1.4){\circle*{0.1}}
\put(3,1.6){\circle*{0.1}}

\end{picture}
\caption{$\textrm{The $2$-vertex-cut $\{ u,x \}$.}$}
\label{D14}
\end{figure}

If $|B_{2}|>1$, then let $x$ denote the neighbour of $v$ in $B_{2}$,
let $\widehat{B_{1}} = B_{1} \cup v$
and let $\widehat{B_{2}} = B_{2} \setminus x$.
Note that $ux$ forms a $2$-vertex-cut where $u$ and $x$ both have just one edge to $\widehat{B_{1}}$
(see Figure~\ref{D14}).
Hence, we can copy case (a) with $B_{1}$ and $B_{2}$ replaced by $\widehat{B_{1}}$ and $\widehat{B_{2}}$, respectively,
to again obtain graphs $B_{1}^{\prime}$ and $B_{2}^{\prime}$ and appropriate augmentations.
It may be that the graph $B_{1}^{\prime}$ will have a $2$-vertex-cut,
so this time we won't be able to use Lemma~\ref{lemma2}
to determine the satisfiability of augmentations of it.
However, we know that we will have $|B_{1}^{\prime}|=4$,
so the number of vertices in any augmentation of $B_{1}^{\prime}$ will be bounded by a constant,
and hence we will be able to determine satisfiability of these augmentations in $O(1)$ time.

\subsection*{Running Time}

We shall now show that the algorithm takes $O \left( |H|^{2.5} \right)$ time.
It is fairly easy to see that the first three stages can be accomplished within this limit
(in fact, they take only $O \left( |H|^{2} \right)$ time),
so we will proceed straight to an examination of Stage 4.

We apply Stage 4 to each of the $2$-vertex-connected blocks derived from Stage~3.
It is easy to see that the total number of vertices in all these blocks is at most $2|H|$,
since each vertex of $H$ will only appear in at most two of these,
so it will actually suffice just to deal with the case when $H$ is itself a $2$-vertex-connected block,
i.e.~when we start Stage 4 with only one $2$-vertex-connected block,
and it has $|H|$ vertices.

During each iteration of Stage 4,
we take a graph $B$
and use it to construct graphs $B_{1}^{\prime}$ and~$B_{2}^{\prime}$,
where $|B_{1}^{\prime}| + |B_{2}^{\prime}| = |B|+2$
and $|B_{2}^{\prime}| < |B|$,
before replacing $B$ with $B_{2}^{\prime}$ and iterating again.
Let $B_{1,1}^{\prime},B_{1,2}^{\prime}, \ldots, B_{1,l}^{\prime}$, for some $l$, 
denote the various graphs that take the role of $B_{1}^{\prime}$ during our algorithm.
Since $|B_{2}^{\prime}| < |B|$,
we can only have at most $|H|$ iterations,
and so we must have $\sum_{i} |B_{1,i}^{\prime}| \leq 3|H|$
(by telescoping, since we always have $|B_{1}^{\prime}| + |B_{2}^{\prime}| = |B|+2$).
We need to apply the algorithm given by Lemma~\ref{lemma2} 
to at most three augmentations of each $B_{1,i}^{\prime}$,
so the total time taken by all such applications will be at most
$3\lambda \sum_{i} \left( |B_{1,i}^{\prime}|^{2.5} \right) \leq 3\lambda \left( \sum_{i}|B_{1,i}^{\prime}| \right)^{2.5} 
= O \left( |H|^{2.5} \right)$.

At the start of each iteration,
we wish to determine whether $B$ has any $2$-vertex-cuts and,
if so, find a minimal one.
Using an algorithm from~\cite{hop} for decomposing a
graph into its so-called `triconnected components',
this takes $O(|B|)=O(|H|)$ time.
It is fairly clear that all other operations involved in an iteration of Stage 4,
aside from applications of Lemma~\ref{lemma2},
can also be accomplished within $O(|H|)$ time,
so (since we recall that there are at most $|H|$ iterations)
this all takes $O \left( |H|^{2} \right)$ time in total
(in fact, by careful bookkeeping, this could be reduced to $O(|H|)$).
Hence,
it follows that the whole algorithm takes $O \left( |H|^{2.5} \right)$ time.

\subsection*{Comments}

By keeping track of all the operations,
the algorithm can be used to find an explicit $4$-regular planar multigraph $G \supset H$ if such a graph exists,
also in $O \left( |H|^{2.5} \right)$ time.
If $H$ is simple,
then we can also obtain a $4$-regular simple planar graph $G^{\prime} \supset H$
without affecting the order of the overall running time,
using the proof of Lemma~\ref{introlemma}.

Aside from looking for a perfect matching during the applications of Lemma~\ref{lemma2},
every part of the algorithm can be accomplished in $O \left( |H|^{2} \right)$ time.
It would therefore be interesting to know if the special structure of the graphs seen in Lemma~\ref{lemma2}
could be exploited to obtain a faster perfect matching algorithm. \\
\\
\\

\end{document}